\documentclass[smallextended,
    envcountsect,]{svjour3}
\smartqed
\usepackage{graphicx}
\journalname{JOTA}

\usepackage{latexsym}
\usepackage{amssymb,amsmath,mathrsfs}
\usepackage{graphicx,float,caption}
\usepackage{tikz-cd}
\usepackage{scalerel} 
\usepackage{color}
\usepackage{graphics}

\usepackage{enumitem}
\setlist{
	leftmargin=*, 
	topsep = 4pt, 
	itemsep = 1pt, 
	widest=viii
}

\usepackage[T1]{fontenc}

\renewcommand{\qed}{\hfill\rule{1ex}{1.5ex}}
\newcommand{\qedex}{\hfill$\triangle$}

\newcommand{\reDeclareMathOperator}[2]{\let#1\undefined \DeclareMathOperator{#1}{#2}}
\newcommand{\reDeclareMathOperatorL}[2]{\let#1\undefined \DeclareMathOperator*{#1}{#2}}

\DeclareMathOperator{\dom}{dom}			
\DeclareMathOperator{\bd}{bd}
\DeclareMathOperator{\gr}{gr}			
\DeclareMathOperator{\cl}{cl}			
\DeclareMathOperator{\wcl}{w-cl}		
\DeclareMathOperator{\inter}{int}		
\DeclareMathOperator{\cone}{cone}		
\DeclareMathOperator{\conv}{conv}		
\DeclareMathOperator*{\argmin}{argmin}	
\DeclareMathOperator*{\argmax}{argmax}	
\DeclareMathOperator*{\Liminf}{Liminf}	

\reDeclareMathOperator{\mod}{mod}
\DeclareMathOperator{\dist}{dist}

\DeclareMathOperator{\interior}{int}

\reDeclareMathOperator{\supp}{supp}

\reDeclareMathOperator{\Proj}{Proj}

\renewcommand{\emptyset}{\varnothing}

\newcommand{\R}{\mathbb{R}}
\newcommand{\N}{\mathbb{N}}

\newcommand{\norm}[1]{\|#1\|}
\newcommand{\abs}[1]{\left|#1\right|}
\newcommand{\product}[1]{\langle #1 \rangle}
\renewcommand{\multimap}{\rightrightarrows}

\newcommand{\xb}{\bar{x}}
\newcommand{\yb}{\bar{y}}

\newcommand{\DA}[1]{\textcolor{black}{#1}}
\newcommand{\blue}[1]{ {\color{black} #1} }

\begin{document}

\title{\blue{Variational and quasi-variational inequalities under local reproducibility: \DA{solution} concept and applications}
}

\titlerunning{Variational and quasi-variational inequalities under local reproducibility}

\subtitle{$\ $\\[5pt]{\small {In honor of Boris Sh. Mordukhovich on the occasion of his of 75th birthday}}}

\author{Didier Aussel and Parin Chaipunya}

\institute{Didier Aussel \at
             Lab. PROMES UPR CNRS 8521\\
             University of Perpignan\\
             Perpignan, France\\
             aussel\@@univ-perp.fr
           \and
             Parin Chaipunya, Corresponding author \at
			 Department of Mathematics\\
             King Mongkut’s University of Technology Thonburi (KMUTT)\\
             Bangkok, Thailand\\
             parin.cha\@@kmutt.ac.th
}

\date{Received: date / Accepted: date}

\maketitle

\begin{abstract}
\blue{Local solutions for variational and quasi-variational inequalities are usually the best type of solutions that could practically be obtained when \DA{in case of lack of} convexity or else when available numerical techniques are too limited for global solutions. \DA{Nevertheless}, the analysis of such problems found in the \DA{literature} seems to be very restricted to the global treatment. Motivated by this fact, in this work, we propose local solution concepts, study their interrelations and relations with global concepts and prove existence results as well as stability of local solution map of parametric variational inequalities. The key ingredient of our results is the new concept of local reproducibility of a set-valued map, which we introduce to explore such local solutions to \DA{quasi-variational inequality} problems. As a by-product, we obtain local solutions to quasi-optimization problems, bilevel quasi-optimization problems and \DA{Single-Leader-Multi-Follower} games.} 
\end{abstract}


\section{Introduction}\label{sec:intro}
\blue{Since} the pioneering works of G. Stampacchia \cite{zbMATH03229929} and Minty \cite{zbMATH03180487}, the concept of variational inequality (VI in short) has been extensively used in the modelling of mechanical systems. But a domain of mathematical analysis in which both of the so-called Stampacchia-type and Minty-type variational inequalities play a central role is the first order analysis of optimization problems. It is indeed well-known that when considering the minimization of a differentiable function over a convex nonempty constraint set, the Stampacchia variational inequality defined by the gradient of the objective function and the constraint set expresses a necessary optimality condition. \blue{The latter becomes necessary and sufficient} when the objective function and the constraint set are convex. This interrelation is still valid when the objective function is not differentiable by replacing the gradient of the objective function by its subdifferential and by considering a set-valued version of Stampacchia variational inequality. The interested reader \DA{can refer} to the book of Facchinei and Pang \cite{Facchinei_Pang_book1,Facchinei_Pang_book2} for more information.

\medskip

More recently the perfect variational reformulation of optimization problems (necessary and sufficient optimality condition) \blue{known} to be limited to convex optimization, has been extended to quasiconvex optimization thanks to the concept of {\em adjusted normal operator} \cite{A_Hadjisavvas05,A_chapter14}. The use of this normal operator allows for example to obtain a sufficient optimality condition for an optimization problem with a continuous quasiconvex function over a nonempty (possibly non convex) constraint set (see eg. Proposition \ref{prop:iff_quasiconvex_optim} below and \cite{A_Hadjisavvas05,A_chapter14}).

\medskip

\blue{We would like to stress that variational and quasi-variational inequalities, widely studied in a global sense, still need more developments \DA{in the local sense} in a similar path to what is done for local minima in optimization problems. A further motivation of this purpose is that the \DA{large majority} of available algorithms consider local minima whenever the objective function misses the convexity property. To our knowledge, only few papers consider local solutions of variational inequality with a remarkable restriction to Minty solutions as in \cite{zbMATH02162400,zbMATH06236761,zbMATH06780576}. \DA{Note that}, the local solutions to Stampacchia variational inequalities and their sensitivity analysis were \DA{defined and} studied in \cite{zbMATH00089169}.} In \cite{zbMATH02125839}, the authors used the concept of local Minty variational inequality as a step to prove the existence of global solutions of Stampacchia variational inequality. This approach has then been extended to equilibrium problem in \cite{zbMATH02162400,Farajzadeh_Zafarani10}.

\medskip

Our aim in this work is to show that local solution concepts for both Minty and Stampacchia variational and quasi-variational inequalities enjoy a rich set of properties. We first investigate the interrelations between the different local concepts and with their global counterparts. Existence of local solutions of variational inequalities are also investigated as well as qualitative stability of the local solutions maps of \blue{parametric} variational inequalities. 

\medskip

When considering local concepts for variational and quasi-variational inequalities, local convexity comes naturally into the hypothesis. But additionally, \blue{inspired by the global treatment in \cite{A_Sagratella17}, we introduce the notion of local reproducible set-valued maps which will play a decisive role in our main results.} Also recall from \cite{A_Sagratella17} that the full solution set of a quasi-variational inequalities for which the constraint set-valued map is reproducible can be determinated by the union of solution sets of the associated variational inequatilies generated by the underlying quasi-variational inequality. A local counterpart of this property is \DA{systematically} under consideration here in this paper.

\medskip

Finally, as applications of the above developments, we give some sufficient conditions for the existence of solution for quasi-optimization and bilevel quasi-variational optimization problems. Let us recall (see e.g.  \cite{Facchinei_Kanzow07,A_Cotrina13,A_Sultana_Vetrivel16}) that quasi-optimization problems are optimization problems in which the constraint set depends on the current value of the variable while bilevel quasi-optimization problem are bilevel optimization problems in which the lower level problem is a quasi-optimization problem. Using a particular case of this last class, we also consider the case of Single-Leader-Multi-Follower games.

\medskip

The paper is organized as follows. In Section \ref{sec:prel}, we set the notation and recall some useful results. In Section \ref{sec:local_def}, the concepts of local solution for variational and quasi-variational are defined and interrelations between the different local solution concepts and their global counterpart are studied \DA{thanks to the use of} the notion of locally reproducible set-valued map. Basic existence results for local solutions of variational inequalities and quasi-variational inequalities are stated in Section \ref{sec:existence} while Section \ref{sec:Stab} is dedicated to qualitative stability of local solution maps of parametrized Stampacchia variational inequalities. These results are then applied in Section \ref{sec:optim} to constrained optimization problems, quasi-optimization problems and finally Single-Leader-Multi-Follower games.

\section{Preliminaries and notations}\label{sec:prel}

Unless otherwise specified, always let $X$ be a Banach space with the norm $\norm{\cdot}$, \blue{$X^{\ast}$ its dual space, and $\langle\cdot,\cdot\rangle$ denotes the duality pair on $X^{\ast} \times X$}. For any $x,y \in X$ the notation $[x,y]$ stands for $\big\{tx+(1-t):t\in[0,1]\big\}$ and the segments $]x,y]$, $[x,y[$, $]x,y[$ are defined by the similar way.

For any given nonempty subset $A \subset X$ and $x \in X$, the distance from $x$ to $A$ is defined by $d(x,A) := \inf \big\{ \|x - y\|: y \in A\big\}$. The notations $\cl(A)$, $\inter(A)$, $\bd(A)$, $\cone(A)$ and $\conv(A)$ stand respectively for the closure, the interior, the boundary, the conical hull and the convex hull of a subset $A$. By $B(x,\rho)$, $\bar{B}(x,\rho)$ and $S(x,\rho)$ respectively, we denote the open ball, the closed ball and the sphere of center $x \in X$ and radius $\rho > 0$.  \blue{We also write 
\[
	\blue{B(A,\rho) := \{x \in X : d(x,A) < \rho\}}
\]
and 
\[
	\blue{\bar{B}(A,\rho) := \{x \in X : d(x,A) \leq \rho\}}
\]
to respectively denote the open and closed balls around the set $A$.} The unit ball and the unit sphere (centered at $0$) on $X$ will be denoted by $\mathbb{B}$ and $\mathbb{S}$, respectively. By a neighbourhood of a point $x \in X$, we mean a set $U \subset X$ that contains an open set that includes $x$.

We next recall the well-known notions concerning set-valued maps. For any set-valued map $T:X \multimap Y$, its domain and the graph are denoted as
\begin{align*}
&\dom T: = \big\{x\in X: T(x)\ne \emptyset\big\},
&\gr T := \big\{(x,y)\in X\times Y: y \in T(x)\big\}.
\end{align*}
The fixed point set of $T$ is defined by $FP(T) := \{x \in X \,|\, x \in T(x)\}$.

\begin{definition}
	Given the maps $ T : X \multimap X^{\ast}$ and $K : X \multimap X$, a point $\bar{x} \in X$ is called:
	\begin{enumerate}[label=-]
		\item a solution to the {\itshape Stampacchia quasi-variational inequality} defined by $T$ and $K$ (briefly $SQVI(T,K)$) if $\bar{x} \in K(\bar{x})$ and there exists $\bar{x}^{\ast} \in T(\bar{x})$ satisfying
		\[
			\langle \bar{x}^{\ast}, y - \bar{x} \rangle \geq 0
		\]
		for all $y \in K(\bar{x})$;
		\item a solution to the {\itshape Minty quasi-variational inequality} defined $T$ and $K$ (briefly $MQVI(T,K)$) if $\bar{x} \in K(\bar{x})$ and for each $y \in K(\bar{x})$ and each $y^{\ast} \in T(y)$, the inequality
		\[
			\langle y^{\ast}, y - \bar{x} \rangle \geq 0
		\]
		holds true.
	\end{enumerate}
	If $K$ is a constant set-valued map with $K(x)$ identical to $C \subset X$ for all $x \in X$, then $SQVI(T,K)$ is reduced to the {\itshape Stampacchia variational inequality} (briefly $SVI(T,C)$). Similarly, the problem $MQVI(T,K)$ reduces to the {\itshape Minty variational inequality} (briefly $MVI(T,C)$).
\end{definition}
When no ambiguity may occurs, we use the same notation for the variational (respectively quasi-variational) inequality and its solution set; for example $SQVI(T,K)$ will both stands for the Stampacchia quasi-variational inequality defined by the couple of maps $(T,K)$ and for the set of its solutions.

\medskip

Recall now some definitions of generalized monotonicity of set-valued maps which will be used thoroughly in the forthcoming. Let $C \subset X$. A set-valued map $T:X\multimap {X^*}$ is said to be 
\begin{enumerate}[label = -]
	\item {\em quasimotonone} on $C$ if, for all $x,y\in C$,
	\[
	\exists~x^*\in T(x)~:~\left\langle x^*,y-x\right\rangle>0
	\implies~\left\langle y^*,y-x\right\rangle\geq0,~~\forall~y^*\in T(y);
	\]
	\item {\em properly quasimonotone} on $C$ if, for every 
	$x_1,x_2,\cdots, x_n\in C$, and all 
	$x\in \conv\left\lbrace x_1,x_2,\cdots,x_n\right\rbrace$, there exists 
	$i\in \left\lbrace 1,2,\cdots,n\right\rbrace$ 
	such that
	\[
	\langle x^*_i,x_i-x\rangle\geq0, \quad\forall x^*_i\in T(x_i);
	\]
	\item {\em pseudomotonone} on $K$ if, for all $x,y\in C$,
	\[
	\exists~x^*\in T(x)~:~\left\langle x^*,y-x\right\rangle\geq0
	\implies~\left\langle y^*,y-x\right\rangle\geq0,\quad \forall~y^*\in T(y).
	\]
\end{enumerate}
It is straightforward to see that pseudomonotonicity implies  properly quasimonotonicity which further implies quasimonotonicity. These concepts of generalized monotonicity provide possibility of the first order analysis of optimization problem with quasiconvex or pseudoconvex objective functions, extending that from the classical convex analysis {(see e.g. \cite{A_Hadjisavvas05,A_chapter14})}.

On the other hand, classical continuity notions for a set-valued operator $T: X\multimap Y$, where $Y$ is a topological vector space, 
are given as follows: $T$ is said to be
\begin{enumerate}[label=-]
	\item {\em lower semicontinuous} at $x_0$ in $\dom T$ if, for any 
	sequence $(x_n)_n$ of $X$ converging to $x_0$ and any element $y_0$ of
	$T(x_0)$, there exists a sequence $(y_n)_n$ of $Y$ converging to $y_0$, with 
	respect to the considered topology on $Y$, and such that $y_n\in T(x_n)$, for any $n$;
	\item {\em upper semicontinuous} at $x_0\in \dom T$ if, for any neighbourhood
	$V$ of $T(x_0)$, there exists a neighbourhood $U$ of $x_0$ such that
	$T(U)\subset V$;
	\item {\em closed} at $x_0\in \dom T$ if, for any sequence 
	$((x_n,y_n))_n \subset \gr T$ converging to $(x_0, y_0)$, one has
	$(x_0, y_0)\in \gr T$;
	\item \blue{{\em closed} if its graph $\gr T$ is a closed set in $X \times X^{\ast}$.}
\end{enumerate}
\blue{It is well-known (e.g. \cite[Proposition 1.4.8]{zbMATH00045282}) that the graph of an upper semicontinuous set-valued map with closed domain and closed values is closed. The converse is true when $T$ is mapped into a compact set.} Whenever $Y$ appearing above is a Banach space equipped with a weak topology, the use of sequences in the above definition should then be replaced by the use of nets (see, e.g., \cite[Proposition 3.2.14]{Run05}), i.e., $T$ is said to be
\begin{enumerate}[label=-]
	\item {\em weakly closed} iff, for any net 
	${(x_\alpha,y_\alpha)}_{\alpha\in \Lambda}\subset \gr T$ such that 
	$(x_\alpha)_\alpha$ weakly converges to $\xb$ and $(y_\alpha)_\alpha$ weakly 
	converges to $\yb$, we have $\yb\in T(\xb)$; 
	\item {\em weakly lower semicontinuous} iff, for any net 
	$(x_\alpha)_{\alpha\in\Lambda}$ weakly converging to $\xb$ and 
	any $\yb\in T(\xb)$, there exists a net $(y_\alpha)_\alpha$
	weakly converging to $\yb$ and such that $y_\alpha\in T(x_\alpha)$, for any $\alpha$.
\end{enumerate}
In the particular case where $Y=X^*$, the weak$^*$ topology should be considered and thus, in both above definitions, the weak convergence of the net 
$(y_\alpha)_\alpha$ is replaced with the $w^*$-convergence.

It turns out that the above notions of continuity for set-valued maps are not always suitable for our use since we will consider cone-valued maps. Hadjisavvas \cite{zbMATH02066991} defined the following two refinements of semicontinuity which is directional in nature and are very well adapted for the study of quasimonotone or pseudomonotone set-valued operators. A set-valued operator $T: X\multimap {X^*}$ is said to be {\em upper sign-continuous} on a convex set $C \subset X$ if, for any $x,y\in C$, the following implication holds:
\begin{equation}\label{eqn:upper_sign_cont}
\forall t\in ]0,1[,~ \inf_{x_t^*\in T(x_t)}\langle x_t^*,y-x\rangle \geq 0
\implies ~ \sup_{x^*\in T(x)}\langle x^*,y-x\rangle \geq 0.
\end{equation}
Note that if $T$ is upper semicontinuous (or even if $T$ is upper hemicontinuous, that is the restriction of $T$ on all lines is upper semicontinuous), then $T$ is upper sign-continuous. 

Finally, \DA{following \cite{A_Hadjisavvas05},} a set-valued map $T: X\multimap {X^\ast}$ is said to be {\em locally upper sign-continuous} on a convex set $C \subset X$ if, for every $x\in C$ there exists a convex neighbourhood $V_x$ of $x$ and an upper sign-continuous operator $S_x : V_x\cap C \multimap X^{\ast}$ with nonempty, convex $w^*$-compact values satisfying $S_x(y)\subset T(y)\setminus\{0\}, \forall\,y\in V_x\cap C$.

\blue{In the recent work of Castellani and Giuli \cite{zbMATH06236761}, a different version of \emph{local} upper sign-continuity \DA{has been} considered without the use of a submap as adopted in the above paragraph. In particular, $T : X \multimap X^{\ast}$ is said to be \emph{locally upper sign-continuous in the sense of Castellani and Giuli} at $x \in X$ if there is a neighbourhood $U_{x}$ of $x$ in which the implication \eqref{eqn:upper_sign_cont} holds for all $y \in C \cap U_{x}$. They also proved that $T$ is locally upper sign-continuous in the sense of Castellani and Giuli at $x$ implies that the same implication \eqref{eqn:upper_sign_cont} holds for all $y \in C$, provided that $C$ is convex. In other words, their local and global upper sign-continuity are equivalent. We remark here the difference between our submap approach (mainly due to the $w^{\ast}$-compactness) to that of Castellani and Giuli \cite{zbMATH06236761} through the following example.
\begin{example}
	Let $X = \R$ and define $T : \R \multimap \R$ by
	\[
		T(x) := \begin{cases}
			]-1,0[	& \text{if $x = 0$} \smallskip\\
			]0,1[	&\text{if $x \neq 0$}.
		\end{cases}
	\]
	Let $x_{0} = 0$. \DA{The map $T$ is clearly (locally) upper sign-continuous at $x_0$ in the sense of Castellani and Giuli since, for any $y<0$ and any $t\in\,]0,1[$ one has $\inf_{x^*_t\in T(x_t)} \langle x^*_t,y-x_0\rangle <0$ and, for any $y>0$ and any $t\in\,]0,1[$, $\sup_{x_{0}^{\ast} \in T(x_{0})} \langle x_{0}^{\ast}, y - x_{0} \rangle = 0$	. Let us now} show that $T$ is not locally upper sign-continuous in the sense of \DA{\cite{A_Hadjisavvas05}}. For \DA{any neighbourhood $U_{x_0}$ of $x_0$ and any submap $S_{x_{0}} : \R \multimap \R$ with $S(y) \subset T(y)\setminus\{0\}$ for all $y \in \R$ and $S(y)$ is nonempty and compact one can find $y\in U_{x_0}$ such that for any $t\in \,]0,1[$, we get
	\[
		\text{$\inf_{x^{\ast}_t \in S(x_t)} \langle x^{\ast}_t, y - x_{0} \rangle > 0$\quad \mbox{ and }\quad $\sup_{x_{0}^{\ast} \in S(x_{0})} \langle x_{0}^{\ast}, y - x_{0} \rangle < 0$.}
	\]
	}
	That is, $T$ contains no compact-valued upper sign-continuous submap. \qedex
\end{example}
}

The above concepts of set-valued maps are strongly motivated from the study of first-order analysis for nonconvex optimization problems. For any function $f : X \to \,]-\infty,+\infty]$ and $x \in X$, the {\em sublevel set} $S_{f} (x)$ and {\em strict sublevel set} $S_{f}^{<}(x)$ of $f$ at $x$ are given respectively by
\[
	S_{f}(x) := \{y \in X \,|\, f(y) \leq f(x) \} \quad\text{and}\quad S_{f}^{<}(x) := \{y \in X \,|\, f(y) < f(x) \}.
\]
For any $x \in X$, we denote $\rho_{x} := \dist(x,S_{f}^{<}(x))$. Following \cite{zbMATH02125839}, the {\em adjusted sublevel set} $S_{f}^{a}(x)$ of $f$ at $x$ is then defined by
\[
	S_{f}^{a}(x) := \begin{cases}
		\bar{B}(S_{f}^{<}(x),\rho_{x}) \cap S_{f}(x)	&\text{if $x \not\in \argmin_{X} f$,} \medskip\\
		S_{f}(x)	&\text{otherwise.}
	\end{cases}
\]
By definition, {a} function $f : X \to \,]-\infty,+\infty]$ is called {\em quasiconvex} if $S_{f}(x)$ is convex for all $x \in X$. However, we also know from \cite{A_chapter14} that $f$ is quasiconvex if and only if $S_{f}^{a}(x)$ is convex for every $x \in X$. Recall that $f$ is called {\em semistrictly quasiconvex} if it is quasiconvex and additionally for any $x,y \in X$, the following condition holds:
\[
	f(y) < f(x)	 \implies [y,x[ \, \subset S_{f}^{<}(x).
\]
It is not difficult to see that for a lower semicontinuous, semistrictly quasiconvex function $f$, we have $\cl(S_{f}^{<}(x)) = S_{f}(x)$ for all $x \not\in \argmin_{X} f$. This fact further implies that $S_{f}^{a}(x) = S_{f}(x)$ for all $x \in X$ for such functions.

Let us recall the notion and properties of an {\em adjusted normal operator}, which will {play} a central tool in our main results.
\begin{definition}[\cite{A_Hadjisavvas05,A_chapter14}]
	To any function $f : X \to \, ]-\infty,+\infty]$, the {\em adjusted normal operator} of $f$ is the set-valued map $N_{f}^{a} : X \multimap X^{\ast}$ defined by
	\[
		N_{f}^{a} (x) := \{x^{\ast} \in X^{\ast} \,|\, \langle x^{\ast}, y - x \rangle \leq 0,\; \forall y \in S_{f}^{a}(x)\}.
	\]
\end{definition}
Note that for each $x \in X$, $N_{f}^{a}(x)$ is a closed convex cone. If $f$ is lower semicontinuous and semistrictly quasiconvex, then $N_{f}^{a}$ can be simplified so that
\[
	N_{f}^{a} (x) = \{x^{\ast} \in X^{\ast} \,|\, \langle x^{\ast}, y - x \rangle \leq 0,\; \forall y \in S_{f}(x)\}.
\]
In this case, we use the notation $N_{f}$ which corresponds to the fact that $S_{f}^{a}$ reduces to $S_{f}$.

{Let us recall some basic properties of the adjusted normal operator.}

\begin{proposition}[\cite{A_Hadjisavvas05,A_chapter14}]\label{prop:prop_normal}
	Let $f : X \to \, ]-\infty,+\infty]$ be any function. Then $N_{f}^{a}$ is quasimonotone. If $f$ is lower semicontinuous and either $f$ is radially continuous or $\dom(f)$ is convex with $\interior(S_{f}^{a}(x)) \neq \emptyset$ for each $x \not\in \argmin_{X} f$, then it holds:
	\begin{enumerate}[label=\upshape\alph*)]
		\item $N_{f}^{a}(x) \setminus \{0\} \neq \emptyset$ on a dense subset of $\dom(f)\setminus\argmin_{X}f$ implies that $f$ is quasiconvex, and
		\item $f$ is quasiconvex implies $N_{f}^{a}(x) \setminus \{0\} \neq \emptyset$ for every $x \in \dom(f)\setminus\argmin_{X}f$.
	\end{enumerate}
\end{proposition}

We shall now state the recent concept of sub-boundarily constant function which turns out to be useful in the study of quasiconvex programs, through the normal approach. Suppose that $f : X \to\, ]-\infty,+\infty]$. We say that $f$ is {\itshape sub-boundarily constant} on $C \subset X$ if for each $x \in C$, the following implication holds:
\[
	f(y) < f(x) \implies [y,x[ \; \cap \interior S_{f}^{a}(x) \neq \emptyset.
\]
Notice that if $f$ is radially continuous and quasiconvex, then $f$ is sub-boundarily constant on $\interior (\dom f )$. If $f$ is semistrictly quasiconvex, then it is sub-boundarily constant if and only if $\interior S_{f(\cdot)} \neq \emptyset$ on $\dom (f) \setminus \argmin f$. It then implies that a continuous semistrictly quasiconvex function is sub-boundarily constant.

\begin{proposition}[\cite{A_CoaVan_Salas21}]\label{prop:Ff}
	If the function $f : X \to \R$ is sub-boundarily constant and quasiconvex, and $F : X \multimap X$ is a set-valued map defined by
	\begin{equation}\label{eqn:Ff_mapping}
	F_{f}(x) := \left\{
	\begin{array}{ll}
	N_{f}^{a}(x) \cap \mathbb{B}     &\text{if $x \in \argmin_{X} f$,} \medskip\\
	\conv ( N_{f}^{a} (x) \cap \mathbb{S} )  &\text{otherwise,}
	\end{array}\right.
	\end{equation}
	then $F_{f}$ is upper sign-continuous and $0 \in F_{f}(x)$ if and only if $x \in \argmin_{X} f$.
\end{proposition}


\section{Local solutions of variational inequalities and quasi-variational inequalities}\label{sec:local_def}

\subsection{Definitions of local solution concepts}\label{subsec:def_local}
Always let $X$ be a Banach space, $T : X \multimap X^{\ast}$ and $K : X \multimap X$ be set-valued maps. We begin this section by precisely stating the definitions of local solutions of several classes of variational inequalities.
\begin{definition}
	A point $\bar{x} \in X$ is called:
\begin{enumerate}[label=-]
	\item a \emph{local solution} to the \emph{Stampacchia quasi-variational inequality} defined by $T$ and $K$ (briefly, $SQVI(T,K)$) if $\bar{x} \in K(\bar{x})$  and there exist a neighbourhood $U$ of $\bar{x}$ and a point $\bar{x}^{\ast} \in T(\bar{x})$ such that 
	\begin{equation}\label{eqn:SVI}
		\langle \bar{x}^{\ast}, y - \bar{x} \rangle \geq 0
	\end{equation}
	for all $y \in K(\bar{x})\cap U$; The set of local solutions of $SQVI(T,K)$ will be denoted by $LSQVI(T,K)$. Now if $K(x)$ is identically $C \subset X$, for all $x \in X$, $LSVI(T,C)$ will stand for the set of local solutions of the Stampacchia variational inequality $SVI(T,C)$;
	\item a \emph{local solution} to the \emph{Minty quasi-variational inequality} defined by $T$ and $K$ (briefly, $MQVI(T,K)$) if $\bar{x} \in K(\bar{x})$ and there exist a neighbourhood $U \ni \bar{x}$ such that for each $y \in K(\bar{x}) \cap U$ and each $y^{\ast} \in T(y)$, the inequality
	\[
	\product{y^{\ast},{y-\bar{x}}} \geq 0
	\]
	holds true; The set of local solutions of $MQVI(T,K)$ will be denoted by $LMQVI(T,K)$. Now if $K(x)$ is identically $C \subset X$, for all $x \in X$, {the set} $LMVI(T,C)$ will stand for the set of local solutions of the Stampacchia variational inequality $MVI(T,C)$.
\end{enumerate}
\end{definition}
It is well known that when the constraint map $K$ is convex valued then there is no need to define a local concept of solution for the Stampacchia quasi-variational inequality $SQVI(T,K)$ since the set of global solutions $SQVI(T,K)$ and the set of local solutions $LSQVI(T,K)$ coincide. But here we aim to consider Stampacchia quasi-variational inequalities in which the constraint sets $K(x)$ are not assumed to be convex but only locally convex. Let us recall that a typical example of locally convex sets which are not convex are the unions of convex sets, a situation which plays a fundamental role for example in disjunctive optimization (see e.g. \cite{AYe08} and references therein). \blue{\DA{In the following example we enlighten} the difference between the local and global solutions of a variational inequality problem (and similarly for a quasi-variational inequality problem).
\begin{example}
	We shall raise a variational inequality from the following linear mixed-integer program
	\[
		\begin{array}{rl}
			\displaystyle\min_{x_{1},x_{2}} & x_{1} + x_{2} \\[.2em]
			\text{s.t.}& \left\{ \begin{array}{l}
				0 \leq x_{1},x_{2} \leq 1 \smallskip\\
				x_{1} \in \mathbb{Z}.
			\end{array}\right.
		\end{array}
	\]
	Obviously, the set of local solution of this problem is $\{(0,0),(1,0)\}$ while the global solution set is $\{(0,0)\}$.

	Using, ahead of time, the Proposition \ref{prop:iff_quasiconvex_optim} below, one obtain an equivalent variational inequality with the map $\DA{T(x) = F_{f}(x) = \{(1,1)\}}$ and constraint set $C = \{0,1\} \times [0,1]$. Then we have
	\[
		LSVI(T,C) = \{(0,0),(1,0)\}
	\]
	while $SVI(T,C)$ consists of only the point \DA{$(0,0)$}.
	
	Further relationships between local solutions of (quasi-)optimization and (quasi-)variational inequality problems are discussed in Section \ref{sec:optim}. \qedex
\end{example}
}

\blue{Note that the local solutions of a Stampacchia variational inequality was already defined in \cite{zbMATH00089169}.} On the other hand the concept of local solution of Minty variational inequality is not new and has been used as an important step in the proof of \cite[Theorem 2.1]{A_Hadjisavvas05} where some sufficient conditions for the existence of global solutions of Minty variational inequality have been derived. 

\subsection{Interrelations between solution concepts}\label{subsec:interrel}
Our aim in this subsection is to precise the interrelations between the different concepts of local solutions and also between the global concepts and the local concepts. To this end, the notion of local reproducibility of a set-valued map which will be introduced below plays a significant role.

Clearly any global solution of the Stampacchia/Minty variational/quasi-variational inequality is also a local solution of the corresponding problem. 

\begin{proposition}\label{prop:links}
Let $T : X\multimap X^*$ be any set-valued map and assume that $K:X\multimap X$ is nonempty and convex valued.
\begin{enumerate}[label=\upshape \alph*)]
	\item If $T$ is pseudomonotone, then $LMQVI(T,K)\subset MQVI(T,K)$;
	\item If $T$ is locally upper sign-continuous on $X$ then 
	\[LMQVI(T,K)\subset LSQVI(T,K).
	\]
\end{enumerate}
\end{proposition}
\begin{proof} It follows the same lines as the one of \cite[Lemma 2.1]{zbMATH02125839} and is included here for sake of completeness. So let $x$ be an element of $LMQVI(T,K)$. Then, there exists a convex neighbourhood $U_x$ of $x$ such that $x\in MQVI(T,K(x)\cap U_x)$.

a) Let us first assume that $T$ is pseudomonotone. Since $K$ is convex-valued, for any $y\in K(x)$, there exists $z=x+t(y-x)$, $t\in]0,1[$ such that $z\in K(x)\cap U_x$. Then, for any $z^*\in T(z)$, $\langle z^*, y-z\rangle=[(1-t)/t]\langle z^*,z-x\rangle\ge 0$, which, by pseudomonotonicity of $T$, yields
$\langle y^*,y-x \rangle\ge 0$, for all $y^*\in T(y)$ implying that $x\in MQVI(T,K)$.

b) Let $y \in K(x)\cap V_x \cap U_x$ where $V_x$ is a convex neighbourhood of $x$ given by the local upper sign-continuity $T$. Suppose that $S_{x}$ is the induced upper-sign continuous submap of $T$ on $K(x)\cap V_x \cap U_x$ from the definition of local upper-sign continuity. The set $K(x) \cap V_x$ is convex and thus, one can find $z\in\,]x,y]$ such that the segment $[z,y]$ is included in $K(x)\cap V_x\cap U_x$ and therefore 
\[\inf_{u\in [z,y]}\inf_{u^*\in S_x(u)} \langle u^*, u-x\rangle \ge 0.
\]
Since the previous inequality holds true for any $y \in K(x)\cap V_x \cap U_x$ and according to the upper sign-continuity of $S_x$ on $V_{x}$, the $w^*$-compactness of $S_x(x)$ and the Sion minimax, 
\[\max_{x^*\in S_x(x)}\inf_{y \in K(x)\cap V_x \cap U_x}\langle x^*y-x\rangle\ge 0
\]
or in other words, $x$ is element of $LSQVI(T,K)=SQVI(T,K)$ since the sets $K(x)$, $U_x$ and $V_x$ are convex. \qed
\end{proof}

As an intermediate step in the analysis between the different concepts of solution, the following diagram summarizes the relationships between global/local solutions of Stampacchia and Minty quasi-variational inequalities.
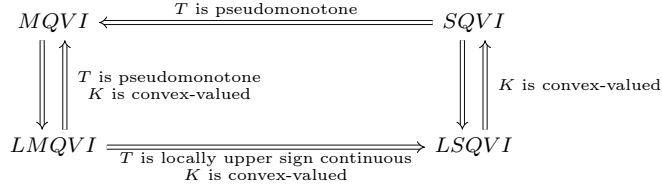
\begin{figure}[H]
	\centering
	\begin{tikzcd}[row sep = 1.2 cm, column sep = 4.2cm]
		MQVI \arrow[r,Leftarrow,"\text{$T$ is pseudomonotone}"] \arrow[d,Rightarrow,xshift=-.5em] \arrow[d,Leftarrow,xshift=.5em,"~\substack{\text{$T$ is pseudomonotone} \\ \text{$K$ is convex-valued}}"] & SQVI \arrow[d,Rightarrow,xshift=-.5em] \arrow[d,Leftarrow,xshift=.5em,"~\text{$K$ is convex-valued}"]   \\
		LMQVI \arrow[r,Rightarrow,"\substack{\text{$T$ is locally upper sign continuous} \\ \text{$K$ is convex-valued}}"'] & LSQVI
	\end{tikzcd}
	\caption{Relationships between concepts of quasi-variational inequalities.}\label{fig:relationships_qvis}
\end{figure}
$\ $\vspace{-1cm}\par

In \cite{A_Sagratella17}, the authors defined the notion of reproducibility of a set-valued map and used it in order to prove that any global solution of a Stampacchia quasi-variational inequality $SQVI(T,K)$ defined with a reproducible constraint map $K$ can be obtained as a solution of a specific Stampacchia variational inequality. That is, in a more formal way,
\begin{equation}\label{eq:A_sagra}
	SQVI(T,K)=\bigcup_{z\in FP(K)} SVI(T,K(z)).
\end{equation}
Let us recall from \cite{A_Sagratella17}, the definition of {(global)} reproducibility of a set-valued map.

\begin{definition}
	Let $C \subset X$. The map $K : X \multimap X$ is said to be \emph{(globally) reproducible} on $C$ if for any $x \in FP(K)\cap C$ and $z \in K(x)$, it holds that $K(x) = K(z)$.
\end{definition}

This reproducility concept could be also used in our local solution context to obtain a ``local counterpart'' of equality \eqref{eq:A_sagra} but we will show that actually a local, and weaker, version of reproducibility will be sufficient to attain this target.


\begin{definition}
	The set-valued map $K:X\multimap X$ is said to be \emph{locally reproducible} at a fixed point $\bar{x} \in FP(K)$ if there exists a neighbourhood $U$ of $\bar{x}$ such that, for any $z \in K(\bar{x}) \cap U$, $K(\bar{x}) \cap U = K(z) \cap U$. Here, $U$ will be called a \emph{neighbourhood of reproducibility at $\bar{x}$}.

	The map $K$ is said to be locally reproducible on a subset $C$ of $X$ if it is locally reproducible at any point $x$ of $C$.
\end{definition}

Let us show, on a simple example, that the local version of reproducibility is clearly weaker than the (global) reproducibility defined in \cite{A_Sagratella17}.
\begin{example}
	Consider the map $K : [-1,1] \multimap [-1,1]$ given by
	\[
	K(x) := \left[ -\sqrt{1 -  x^{2}}, \sqrt{1 - x^{2}} \right]
	\]
	for $x \in [-1,1]$. Then $FP(K) = \left[ -\frac{1}{\sqrt{2}},\frac{1}{\sqrt{2}} \right]$ and it is not difficult to see that $K$ is not reproducible. Indeed, observe that $0 \in K\left(\frac{1}{\sqrt{2}}\right) = \left[ -\frac{1}{\sqrt{2}}, \frac{1}{\sqrt{2}} \right]$ but $K(0) = [-1,1] \neq K\left(\frac{1}{\sqrt{2}}\right)$. However $K$ is clearly locally reproducible on $\left] -\frac{1}{\sqrt{2}}, \frac{1}{\sqrt{2}} \right[$. \qedex
\end{example}

Another possible definition for the local reproducibility of a set-valued map could have been the following
	\begin{quote}
		The map $K$ is said to be \emph{locally reproducible} (in the submap sense) at a fixed point $\bar{x} \in FP(K)$ if there exist a neighbourhood $U$ of $\bar{x}$ and a submap $K_{\bar{x}}$ of $K$ for which $K_{\bar{x}}$ is reproducible on $U$.
	\end{quote}
The interrelation between these two definitions is given in the proposition below.
\begin{proposition}
	If $K$ is locally reproducible at a fixed point $\bar{x} \in FP(K)$, then $K$ is locally reproducible in the submap sense at $\bar{x}$. 
\end{proposition}
\begin{proof}
	Let $U$ be a neighbourhood of reproducibility at $\bar{x}$. Define $K_{\bar{x}} : X \multimap X$ by
	\[
	K_{\bar{x}}(z) := \left\{
	\begin{array}{ll}
		K(z)        &\text{if $z \not\in U$,} \\[.5em]
		K(\bar{x}) \cap K(z)\cap U  &\text{if $z \in U$.}
	\end{array}
	\right.
	\]
	Taking any $y \in FP(K_{\bar{x}}) \cap U$, we get $y \in K(\bar{x}) \cap U$ by the above definition of $K_{\bar{x}}$. Consequently, the local reproducibility at $\bar{x}$ implies that $K(y) \cap U = K(\bar{x}) \cap U$. Again, by the local reproducibility at $\bar{x}$, the following holds for any $z \in K_{\bar{x}}(y)$:
	\[
	K_{\bar{x}}(y) = K(\bar{x}) \cap K(y) \cap U = K(\bar{x}) \cap K(\bar{x}) \cap U = K(\bar{x}) \cap K(z) \cap U = K_{\bar{x}}(z)
	\]
	thus ensuring the reproducibility of $K_{\bar{x}}$ on $U$, as required. \qed
\end{proof}

The converse of the above proposition is not generally true, as is depicted in the following example.
\begin{example}
	Consider a set-valued map $K : \R \multimap \R$ defined by
	\[
	K(x) := \left\{
	\begin{array}{ll}
		[0,1]   &\text{if $x < 1$,} \\[.5em]
		[0,2]   &\text{if $x \geq 2$.}
	\end{array}\right.
	\]
	and the fixed point $\bar{x} := 1 \in FP(K)$. Clearly, $K$ is not locally reproducible at $\bar{x}$. Now for any $\varepsilon \in (0,1)$, we take a neighbourhood $U_{\varepsilon} = [1-\varepsilon,1+\varepsilon]$ of $\bar{x}$ and define a map $K_{\bar{x}}^{\varepsilon} : U_{\varepsilon} \multimap X$ by $K_{\bar{x}}^{\varepsilon} (u) := \left[1-\varepsilon,1\right]$ for any $u \in U_{\varepsilon}$. Then $K_{\bar{x}}^{\varepsilon}$ is a reproducible submap of $K$ on $U_{\varepsilon}$. 
	This shows that $K$ contains a reproducible submap on any neighbourhood of $\bar{x}$ (see Figure \ref{fig:submap_repro_not_loc_repro}), while not being locally reproducible at $\bar{x}$. Also notice that $K$ is upper semicontinuous, hence providing an example of upper semicontinuous set-valued map which is not locally reproducible. \qedex
	\begin{figure}[H]
		\centering
		\includegraphics[width = .65\textwidth]{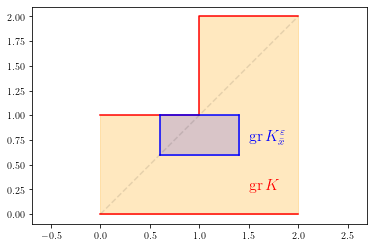}
		\caption{The map $K$ contains a reproducible submap on any neighbourhood of $\bar{x}$ but is not locally reproducible at $\bar{x}$.}
		\label{fig:submap_repro_not_loc_repro}
	\end{figure}
\end{example}

{As proved in the proposition below, and in contrast to what occurs for the global reproducibility, every set-valued map is, roughly speaking, locally reproducible at any fixed point lying in the interior of its graph.}
\begin{proposition}\label{prop:loc.repro_int.gr}
	Let $\bar{x} \in FP(K)$ be such $(\bar{x},\bar{x}) \in \interior (\gr K)$. Then $K$ is locally reproducible at $\bar{x}$.
\end{proposition}
\begin{proof}
	Since $(\bar{x},\bar{x}) \in \interior (\gr K)$, there exists $r > 0$ such that the ball
	\[
	B_{X\times X} ((\bar{x},\bar{x}),r) := \left\{(u,v) \in X \times X \,|\, \norm{u-\bar{x}}^{2} + \norm{v-\bar{x}}^{2} < r^{2} \right\}
	\]
	is contained in $\gr K$. This would mean 
	\begin{equation}\label{eqn:ball_inclusion}
		B_{X\times X}^{\infty} \left((\bar{x},\bar{x}),\frac{r}{\sqrt{2}} \right) \subset B_{X\times X} ((\bar{x},\bar{x}),r) \subset \gr K,
	\end{equation}
	where one may recall that
	\[
	B_{X\times X}^{\infty} \left((\bar{x},\bar{x}),\frac{r}{\sqrt{2}} \right) := \left\{ (u,v) \in X \times X \,\left|\, \max \{\norm{u-\bar{x}},\norm{v-\bar{x}}\} < \frac{r}{\sqrt{2}} \right.\right\}.  
	\]
	Let $U := B_{X} \left(\bar{x},\frac{r}{\sqrt{2}}\right)$ and due to \eqref{eqn:ball_inclusion}, we get $U \subset K(\bar{x})$ which further gives $K(\bar{x})\cap U = U$. Next, take any $y \in K(\bar{x}) \cap U = U$ and any $v \in U$. One has $\max\{ \norm{y-\bar{x}},\norm{v-\bar{x}}\} < \frac{r}{\sqrt{2}}$, showing that $(y,v)$ belongs to $B_{X\times X}^{\infty} \left( (\bar{x},\bar{x}) , \frac{r}{\sqrt{2}} \right)$ and hence to $\gr K$. Thus $U \subset K(y)$ so that $K(y) \cap U = U = K(\bar{x}) \cap U$, completing the proof. \qed
\end{proof}
	
The class of globally reproducible maps of \cite{A_Sagratella17} are quite specific due to the technical hypotheses which were used in those results (see e.g. Propositions 2 and 3). On the contrary, due to the above proposition, we can prove that the local reproducibility can be obtained for a very large (and important) class of set-valued maps, namely the ones defined by {\em separable inequality constraint}. 

\begin{corollary}\label{cor:SIC}
	Let $g,h : X \to \R$ be two functions and $\bar{x} \in X$ be such that $g(\bar{x}) < h(\bar{x})$, $g$ is upper semicontinuous at $\bar{x}$ and $h$ is lower semicontinuous at $\bar{x}$. Then the map $K : X \multimap X$ defined by
	\[
	K(x) := \{y \in X \,|\, g(y) \leq h(x) \}, \qquad \text{for $x \in X$,} 
	\]
	is locally reproducible at $\bar{x}$.
\end{corollary}
\begin{proof}
	The point $\bar{x}$ is clearly a fixed point of $K$. Take any $\varepsilon \in\, ]0,(h(\bar{x})-g(\bar{x}))/2[$. By the upper semicontinuity of $g$ at $\bar{x}$, there exists a neighbourhood $U$ of $\bar{x}$ for which $g(y) < h(\bar{x}) - \varepsilon$, for all $y \in U$. On the other hand, the lower continuity of $h$ at $\bar{x}$ implies the existence of another neighbourhood $V$ of $\bar{x}$ such that $h(z) > h(\bar{x}) - \varepsilon$, for each $z \in V$. Altogether, this means $U \subset K(z)$, for all $z \in V$, and hence $(\bar{x},\bar{x}) \in V \times U \subset \gr K$. By Proposition \ref{prop:loc.repro_int.gr}, $K$ is locally reproducible at $\bar{x}$. \qed
\end{proof}

The remaining of this section would be dedicated to the consequence of local reproducibility, which renders the local solution of  quasi-variational inequalities throught the ordinary variational inequalitiy. More precisely, it is clear that
\begin{equation}\label{eq:union_LSQVI}
	LSQVI(T,K) \subset \bigcup_{z \in FP(K)} LSVI(T,K(z)\cap U_{z})
\end{equation}
and also
\begin{equation}\label{eq_union_LMQVI}
	LMQVI(T,K) \subset \bigcup_{z \in FP(K)} LMVI(T,K(z)\cap U_{z}),
\end{equation}
where $U_{z}$ denotes the neighbourhood of reproducibility at $z \in FP(K)$. It is not difficult to see that the reverse inclusions are not satisfied in general (see e.g. \cite{A_Sagratella17}). The aim of the propositions below are to show that the local reproducibility of the constraint map $K$ allows the reverse inclusion in \eqref{eq:union_LSQVI} and \eqref{eq_union_LMQVI}.

\begin{proposition}\label{prop:LSVI_implies_LSQVI}
	Suppose that $K$ is locally reproducible at a fixed point $z \in FP(K)$ with an open neighbourhood of reproducibility $U_{z}$. Then each local solution $\bar{x}$ of $SVI(T,K(z) \cap U_{z})$ is a local solution of $SQVI(T,K)$. {As a consequence,
	\[
		LSQVI(T,K) = \bigcup_{z \in FP(K)} LSVI(T,K(z)\cap U_{z})
	\]
}
\end{proposition}
\begin{proof}
	Since $\bar{x}$ is a local solution of $SVI(T,K(z) \cap U_{z})$, it belongs to $K(z) \cap U_{z}$. The local reproducibility implies that $K(z) \cap U_{z} = K(\bar{x}) \cap U_{z}$. Since $U_{z}$ is open, it is also a neighbourhood of $\bar{x}$ and we consequently have $\bar{x} \in K(\bar{x})$ and there exist $\bar{x}^{\ast} \in T(\bar{x})$ such that $\langle \bar{x}^{\ast}, y - \bar{x} \rangle \geq 0$ for all $y \in K(\bar{x})\cap U_{z}$. \qed
\end{proof}

\begin{proposition}\label{prop:rep_VI_in_QVI}
	Suppose that $K$ is locally reproducible at a fixed point $z \in FP(K)$ whose neighbourhood of reproducibility at $z$ contains a local solution $\bar{x}$ of $MVI(T,K(z))$ for some $T : X \multimap X^{\ast}$. Then $\bar{x}$ is a local solution of $MQVI(T,K)$. {As a consequence,
	\[
		LMQVI(T,K) = \bigcup_{z \in FP(K)} LMVI(T,K(z)\cap U_{z}),
	\]
}
\end{proposition}
\begin{proof}
	Suppose that $U$ is the neighbourhood of reproducibility at $z$ containing a local solution $\bar{x}$ of $MVI(T,K(z))$. Let $U' \subset U$ be a neighbourhood of $\bar{x}$ such that for all $y \in K(z)\cap U'$ and $y^{\ast} \in T(y)$, we have $\langle y^{\ast},y-x \rangle \geq 0$. Since $\bar{x} \in K(z)\cap U$, the local reproducibility of $K$ at $z$ implies $K(\bar{x}) \cap U = K(z) \cap U$. Since $U' \subset U$, we obtain $K(\bar{x}) \cap U' = K(z) \cap U'$, showing that $\bar{x}$ locally solves $MQVI(T,K)$. \qed
\end{proof}

\begin{corollary}\label{cor:rep:MVI_in_MQVI}
	Suppose that $K$ is locally reproducible at a fixed point $z \in FP(K)$ with the neighbourhood of reproducibility $U$. If $MVI(T,K(z) \cap U)$ has a solution, then $MQVI(T,K)$ has a local solution.
\end{corollary}

\section{Existence results}\label{sec:existence}

In this section, we present existence theorems for local solutions of both SVIs and MVIs under locally convex assumption on the constraint set and quasimonotone-type assumptions on the operator. The results for {both} VIs are then extended to QVIs exploiting the {concept} of local reproducibility.

{Let us first prove} an existence result for a local solution of a {Stampacchia variational inequality}.
\begin{theorem}\label{thm:existence_LSVI}
	Let $C$ be a nonempty, locally convex subset of a reflexive Banach space $X$ and $T : X \multimap X^{\ast}$ a locally upper sign-continuous quasimonotone operator. Then $LSVI(T,C) \neq \emptyset$.
\end{theorem}
\begin{proof}
	Take any $x \in C$ and let $U \subset X$ be a neighbourhood of $x$ such that $C \cap U$ is convex and bounded. Together with the reflexivity of $X$, all the assumptions \cite[Theorem 2.1]{zbMATH02125839} {(actually a \emph{compact} version of it)} are verified so that we have $SVI(T,C \cap U) \neq \emptyset$. This is equivalent to saying $LSVI(T,C) \neq \emptyset$. \qed
\end{proof}

The existence of a local solution of a MVI was proved in \cite{zbMATH02125839}, stating under the convexity of $C$ and the quasimonotonicity of $T$ that $LMVI(T,C) \neq \emptyset$ provided that $T$ is not properly quasimonotone. Otherwise, if $C$ is closed and convex, $T$ is properly quasimonotone and a coercivity condition is satisfied, then the global solution in $MVI(T,C)$ exists (see \cite{zbMATH01367774}).

Here, we show that the closedness of $C$ in the above statement is superfluous and the nonemptiness of $MVI(T,C)$ can be proved alternatively using the transfer method of \cite{zbMATH00090657}. This weaker assumption also enables us to show the existence of a local solution in $LMVI(T,C)$.

Let us recall the notion of transfer closedness and its equivalent formulation from \cite{zbMATH00090657}. 
\begin{definition}[\cite{zbMATH00090657}]
	A set-valued map $G$ from a topological space $X$ into another topological space $Y$ is said to be {\itshape transfer closed-valued} if for any $x \in X$ and $y \not\in G(x)$, there exists $x' \in X$ such that $y \not\in \cl G(x')$.
\end{definition}
\begin{lemma}[\cite{zbMATH00090657}]
	A set-valued map $G$ from a topological space $X$ into another topological space $Y$ is transfer closed-valued if and only if $\bigcap_{x \in X} \cl G(x) = \bigcap_{x \in X} G(x)$.
\end{lemma}
The following lemma is required in the proof of the existence theorem for global and local solutions of a MVI in Theorem \ref{thm:MVI_existence} below. 
\begin{lemma}\label{lem:transfer_closed_criteria}
	Let $X$ and $Y$ be two topological spaces, $h : X \times Y \to \R$ be a bifunction satisfying $h(x,\cdot)$ is continuous and $h(x,x) = 0$ for each $x \in X$.
	Then the map $G : X \multimap Y$ given by
	\[
		G(x) := \{y \in Y \,|\, h(x,y) \leq 0\}
	\]
	is transfer closed-valued.
\end{lemma}
{It's proof, which derived directly from definition, is included for sake of completeness.}
\begin{proof}
	First, note that $G(x) \neq \emptyset$ for all $x \in X$. Suppose now that $x \in X$ and $y \not\in G(x)$, which means $h(x,y) > 0$. This clearly implies that $y$ is not an accumulation point of $G(x)$, so that $y \not\in \cl G(x)$. \qed
\end{proof}

\begin{theorem}\label{thm:MVI_existence}
	Let $C$ be a nonempty subset of a Banach space $X$ and $T : X \multimap X^{\ast}$ be quasimonotone. 
	\begin{enumerate}[label = \upshape\alph*)]
		\item\label{itm:MVI_existence} If $C$ is convex, $T$ is properly quasimonotone, and the following coercivity condition holds
		\begin{equation}\label{eqn:coercivity_for_MVI}
			\begin{array}{c}
			\text{there exist a weakly compact set $W \subset X$ and $x_{0} \in W$ such that}\\
			\forall x \in C \setminus W,\, \exists x_{0}^{\ast} \in T(x_{0}):\, \langle x_{0}^{\ast}, x - x_{0} \rangle > 0,
			\end{array}
		\end{equation}
		then $MVI(T,C) \neq \emptyset$.

		\item\label{itm:LMVI_existence} If $C$ is locally convex and $X$ is reflexive, then $LMVI(T,C) \neq \emptyset$.
	\end{enumerate}
\end{theorem}
\begin{proof}
	Let us first prove \ref{itm:MVI_existence}. Define a set-valued map $G : C \multimap C$ by
	\[
		G(x) := \{ y \in C \,|\, \langle x^{\ast}, y - x \rangle \leq 0,\, \forall x^{\ast} \in T(x) \}.
	\]
	It is clear that {$G$} has nonempty values{, since one always has $x\in G(x)$, even if $T(x)=\emptyset$}. Due to the proper quasimonotonicity of $T$, {for each} $x_{1},\cdots,x_{n} \in C$ and $y \in \conv\{x_{1},\cdots,x_{n}\}$, it holds $y \in \bigcup_{i=1}^{n} G(x_{i})$. This means {that the set-valued map} $\wcl G$, defined by $x \in C \mapsto \wcl G(x)$, is a weakly closed-valued KKM map. Moreover, the coercivity assumption shows that $\wcl G(x_{0}) \subset W$, implying that it is weakly compact. The Ky Fan's Theorem (see \cite{zbMATH03152769}) then implies that $\bigcap_{x \in X} \wcl G(x) \neq \emptyset$.

	Now, we put $h(x,y) := \sup_{x^{\ast} \in T(x)} \langle x^{\ast}, y - x \rangle$ for any $x,y \in C$. It is easy to see that $h(x,x) = 0$ and $h(x,\cdot)$ is weakly continuous. Moreover, $G(x)$ has the following alternative expression
	\[
		G(x) = \{y \in C \,|\, h(x,y) \leq 0 \}.
	\]
	Using Lemma \ref{lem:transfer_closed_criteria}, the map $G$ is transfer weakly closed-valued. Finally, we obtain that $\bigcap_{x \in X} G(x) = \bigcap_{x \in X} \wcl G(x)$ is nonempty {but} its elements are solutions of $MVI(T,C)$. We have arrived at the desired conclusion.

	Next we prove \ref{itm:LMVI_existence}. Note that if $T$ is not properly quasimonotone, then the nonemptiness of $LMVI(T,C)$ follows from \cite[Proposition 2.1]{zbMATH02125839}. It remains to prove the existence in the case where $T$ is properly quasimonotone. Let $x \in C$ and $U$ be a bounded open neighbourhood of $x$ in which $C \cap U$ is convex. Since $C \cap U$ is bounded, we may choose a weakly compact set $W$ large enough so that $(C \cap U) \setminus W$ is empty. Hence the coercivity condition \eqref{eqn:coercivity_for_MVI} automatically holds true for $C \cap U$. The existence of a local solution in $LMVI(T,C)$ follows by applying part \ref{itm:MVI_existence} to the problem $MVI(T,C \cap U)$. \qed
\end{proof}

Finally, thanks to the local reproducibility condition, we can easily extend the above existence results to quasi-variational inequalities.
\begin{theorem}
	Let $X$ be a reflexive Banach space, $K : X \multimap X$ be a locally convex-valued, locally reproducible map on $FP(K)$ with $FP(K) \neq \emptyset$, and $T : X \multimap X^{\ast}$ a quasimonotone operator.
	Then $LMQVI(T,K) \neq \emptyset$. Moreover, if $T$ is locally upper sign-continuous, then $LSQVI(T,K) \neq \emptyset$.
\end{theorem}
\begin{proof}
	Take any $z \in FP(K)$ and an open neighbourhood of reproducibility $U_{z}$ at $z$. Then $K(z) \cap U_{z}$ {can be assumed, wlog, to be} locally convex and hence $LMVI(T,K(z)\cap U_{z}) \neq \emptyset$ by Theorem \ref{thm:MVI_existence}. {Then the $LMQVI(T,K)\neq\emptyset$ thanks to Proposition \ref{prop:LSVI_implies_LSQVI}.} The nonemptiness of $LSQVI(T,K) \neq \emptyset$ under upper sign-continuity condition can be proved similarly {but using Theorem \ref{thm:existence_LSVI} and Proposition \ref{prop:rep_VI_in_QVI}.} \qed
\end{proof}
\begin{remark}
	Actually we can draw a stronger conclusion from the above proof. Indeed, for any $z\in FP(K)$, we obtain at least one local solution of $MQVI(T,K)$ in its neighbourhood.
\end{remark}

\section{Stability of local solutions}\label{sec:Stab}

Our aim in this section is to prove some stability properties forlocal solutions to Stampacchia variational inequalities subject to perturbations. 
These stability results will play a central role in Section \ref{sec:optim} to prove the existence of solutions of a specific Single-Leader-Multi-Follower game $(SLMFG)$, see Theorem \ref{thm:exist_SLMFG}.

Many studies has been done in the litterature and both qualitative (closedness, semi-continuity, etc.) and quantitative (H\"older/Lipschitz property) have been proved for such variational inequality, with perturbations on the map $T$, on the constraint map $K$ or on both (see e.g. \cite{Ait-Mansour_A06,zbMATH05351708,Facchinei_Pang_book1} and references therein). \blue{In particular, the stability in the sense of differentiability of local solutions of a variational inequality problem has been studied by \cite{zbMATH00089169}. However, the structure and approach of \cite{zbMATH00089169} are far different from what we study in this present paper.}

\blue{Let us now fix the framework of our analysis.} Suppose that $\Lambda$ and $U$ are Banach spaces of perturbation parameters and consider the perturbed set-valued maps $T:X\times\Lambda \multimap X^*$ and $K:U\multimap X$. The local solution maps of the perturbed variational inequalities are then defined as follows:
\begin{itemize}
	\item [- ] the {\em local solution map} $LSVI^*_{pert}$ defined, for any $(\lambda,\mu)\in \Lambda\times U$, by $LSVI^*_{pert}(\lambda,\mu)$ being the set of $\bar{x}\in K(\mu)$ for which there exist a neighbourhood $U_{\bar{x}}$ of $\bar{x}$ and $\bar{x}^{\ast}\in T(\bar{x},\lambda)$ with $\bar{x}^{\ast}\neq 0$ such that 
	\[\langle \bar{x}^*, y-\bar{x}\rangle \ge 0,\quad\forall\,y\in K(\mu)\cap U_{\bar{x}};
	\]
	\item [- ] the {\em weak-int-local solution map} $LSVI_{pert}^{w,int}$ defined, for any $(\lambda,\mu)\in \Lambda\times U$, by $LSVI^{w,int}_{pert}(\lambda,\mu)$ being the set of $\bar{x} \in K(\mu)$ such that there exists a neighbourhood $U_{\bar{x}}$ of $\bar{x}$ for which, for any $y \in \interior K(\mu) \cap U_{\bar{x}}$, there exists $\bar{x}^{\ast} \in T(\bar{x},\lambda)\setminus\{0\}$ with $\langle \bar{x}^{\ast}, y - \bar{x} \rangle > 0$.
\end{itemize} 
Note that, as commonly used in the literature on variational inequalities, the ``weak'' terminology stands for weak solutions of the Stampacchia variational inequality, that is solution $\bar{x}$ for which the associated operator $\bar{x}^*$ could depend on the vector $y$. Moreover the notation ``$*$'' means that the considered $\bar{x}^*\in T(\bar{x},\lambda)$ cannot be null.

Let us now give the definition of a ``regularity-type'' property, inspired from \cite{Ait-Mansour_A06} and \cite{A_CoaVan_Salas21}, on the couple of  maps $(T,K)$ defining the perturbed variational inequalities considered in this section.
\begin{definition}\label{def:dually}
	The couple $(T,K)$ of set-valued maps, with $T : X \times \Lambda \multimap X^{\ast}$ and $K : U \to X$ will be said to be {\em int-dually lower semicontinuous} on $X^2\times\Lambda\times U$ if for all sequence $\left(y_n,z_n,\lambda_n,\mu_n\right)_n$ converging to $(y,z,\lambda,\mu)\in X^2\times\Lambda\times U$ with $y \in \interior K(\mu)$, $z \in K(\mu)$ and, for any $n$, $y_{n} \in \interior K(\mu_{n})$,  $z_{n} \in K(\mu_{n})$, it holds
	\[
	\sup_{y^{\ast} \in T(y,\lambda)\setminus\{0\}} \langle y^{\ast}, z - y \rangle \leq \liminf_{n} \sup_{y^{\ast}_{n} \in T(y_{n},\lambda_{n})\setminus\{0\}} \langle y^{\ast}_{n}, z_{n} - y_{n} \rangle.
	\]	
\end{definition}
As observed in \cite{A_CoaVan_Salas21}, considering regularity conditions on the couple $(T,K)$ instead of on each map $T$ and $K$ separately, as classically done, allows to weaken the hypotheses used to prove some stability results for the solution maps. 

\begin{theorem}\label{thm:beta}
	Using the above notations, suppose that for any $(\lambda,\mu) \in \Lambda \times U$, the set $LSVI^{w,int}_{pert}(\lambda,\mu)$ is nonempty and the following conditions are fulfilled:
	\begin{enumerate}[label=\upshape\roman*)]
		\item for all $\mu \in U$, $K(\mu)$ is locally convex with nonempty interior;
		\item\label{assmp:loc_upper_sign} for all $(\lambda,\mu) \in \Lambda \times U$, $T(\cdot,\lambda)$ is quasimonotone on $K(\mu)$ and, for any $x \in \dom T(\cdot,\lambda)$ there exists a neighbourhood $V_{x}$ of $x$ and an upper sign-continuous map $\Phi_{x}(\cdot,\lambda) : V_{x} \cap K(\mu) \multimap X^{\ast}$ with nonempty $w^{\ast}$-compact values satisfying $\Phi_{x}(y,\lambda) \subset T(y,\lambda)\setminus\{0\}$ for all $y \in V_{x} \cap K(\mu)$.
		\item\label{hyp:lsc-like} the couple $(T,K)$ is int-dually lower semicontinuous on
		$X^2\times \Lambda\times U$; 
		\item for any sequence $(\mu_{n})_n$ converging to $\mu$, the sequence $(K(\mu_{n}))_n$ Mosco converges to $K(\mu)$.
	\end{enumerate}
	Then, the map $LSVI^{w,int}_{pert}$ is closed. If, additionally, in assumption \ref{assmp:loc_upper_sign}, the operators $\Phi_{x}(\cdot,\lambda)$ are convex valued (that is the operator $T$ is locally upper sign-continuous on $K(\mu)$), then $LSVI^*_{pert}$ is closed.
\end{theorem}

\blue{We would like to make the following remarks for the above result with respect to the global settings in the close literature before formulating its proof.
\begin{enumerate}[label=\alph*)]
	\item\label{rmk:itm_a} In the above theorem, if $K(\mu)$ is convex (instead of locally convex), then one could refer to Figure \ref{fig:relationships_qvis} and see that local solutions in $LSVI^{\ast}_{pert}$ becomes global.
	\item If the situation in \ref{rmk:itm_a} occurs, Theorem \ref{thm:beta} reduces to a result of Ait Mansour and Aussel \cite[Theorem 4.2]{zbMATH05351708}.
	\item \blue{A similar global results were deduced by Aussel and Cotrina \cite{zbMATH05900884,A_Cotrina13} under different assumptions. For instance, $K$ is assumed to by closed and lower semi-continuous in \cite[Proposition 3.1]{zbMATH05900884} and $T$ satisfies a different kind of continuity in \cite[Proposition 4.3]{A_Cotrina13}.}
\end{enumerate}
}

Now, to prove Theorem \ref{thm:beta}, the following lemma is a prerequisite.

\begin{lemma}\label{lem:Mosco_int_int}
	Suppose that $(S_{n})_n$ is a sequence of locally convex subsets of $X$ with nonempty interiors and that $S$ is a subset of $X$ such that $S \subset \Liminf_{n} S_{n}$.
	Then $\interior S \subset \Liminf_{n} \interior S_{n}$.
\end{lemma}
\begin{proof}
	Without loss of generality one can assume that $\interior S\neq \emptyset$. Suppose that $z \in \interior S \subset \Liminf_{n} S_{n}$ and let $(z_{n})_n$ be a sequence converging to $z$ with $z_{n} \in S_{n}$ for each $n \in \N$. Actually, we have $z_{n} \in \interior S$ for sufficiently large $n \in \N$. In other words, $\interior S$ intersects with every open neighbourhood of $z_{n}$ whenever $n > N$, for some $N > 0$. By the local convexity of each $S_{n}$, we may choose a positive sequence $(\rho_{n})_n$ converging to $0$ and such that $B(z_{n},\rho_{n}) \cap S_{n}$ is convex for every $n > N$. Taking any $z'_{n} \in B(z_{n},\rho_{n}) \cap S_{n}$ for each $n > N$, we obtain that $z'_{n}$ converges to $z$. We therefore get
	\begin{align*}
		z &\in \Liminf_{n} (B(z_{n},\rho_{n}) \cap S_{n}) = \Liminf_{n} \cl(B(z_{n},\rho_{n})\cap S_{n})\\
		& =  \Liminf_{n} \cl\interior (B(z_{n},\rho_{n})\cap S_{n}) = \Liminf_{n} \interior (B(z_{n},\rho_{n})\cap S_{n}) \\
		& = \Liminf_{n} (B(z_{n},\rho_{n}) \cap \interior S_{n})  \subset \Liminf_{n} \interior S_{n}
	\end{align*}
	thus completing the proof.
	\qed
\end{proof}
\begin{proof}{\em of Theorem \ref{thm:beta}}
	Let $(\lambda_{n},\mu_{n})_n$ and $(x_{n})_n$ be two sequences in $\Lambda \times U$ and $X$ respectively, and converging respectively to $x$ and $(\lambda_{0},\mu_{0})$ with
	\[
	x_{n} \in LSVI^{w,int}_{pert}(\lambda_{n},\mu_{n})\mbox{, for all } n \in \N.
	\]
	Let $V_{x}$ be a neighbourhood of $x$ such that $V_{x} \cap K(\mu_{0})$ is convex and $\Phi_{x}(\cdot,\lambda_{0}) : V_{x} \cap K(\mu_{0}) \multimap X^{\ast}$ be a submap of $T(\cdot,\lambda_{0})$ which is upper sign-continuous and has nonempty $w^{\ast}$-compact values.
	
	By the Mosco convergence of the sequence $(K(\mu_{n}))_n$ and the fact that $x_{n} \in K(\mu_{n})$ for each $n \in \N$, one gets $x \in K(\mu_{0})$. Take any $y \in V_{x} \cap K(\mu_{0}) \setminus \{x\}$. The convexity of $V_{x} \cap K(\mu_{0})$ ensures that it contains the segment $[y,x[$. For $t \in ]0,1]$, let $z_{t} := ty + (1-t)x$. By the Mosco convergence and Lemma \ref{lem:Mosco_int_int}, there is a sequence $(z_{n})_n$ converging to $z_{t}$ with $z_{n} \in \interior K(\mu_{n})$, for any $n \in \N$. Since, for each $n \in \N$, $x_{n} \in LSVI^{w,int}_{pert}(\lambda_n,\mu_n)$, we may find $x^{\ast}_{n} \in T(x_{n},\lambda_{n}) \setminus \{0\}$ satisfying
	\[
	\langle x^{\ast}_{n}, z_{n} - x_{n} \rangle \geq 0.
	\]
	Since $z_{n} \in \interior K(\mu_{n})$ and $x_n^*\neq 0$, we may assume, without loss of generality, that the above inequality holds strictly for all $n \in \N$. Combining with the quasimonotonicity of $T(\cdot,\lambda_{n})$, one obtains
	\[
	\langle z^{\ast}_{n}, z_{n} - x_{n} \rangle \geq 0, \quad \forall z^{\ast}_{n} \in T(z_{n},\lambda_{n})\setminus\{0\}.
	\]
	Letting $n \to \infty$ and using the hypothesis \ref{hyp:lsc-like}, we have
	\[
	\langle z^{\ast}_{t}, x - z_{t} \rangle \leq 0
	\]
	for any $z^{\ast}_{t} \in T(z_{t},\lambda_{0})\setminus\{0\}$. For any $t \in\, ]0,1]$ and any $z^{\ast}_{t} \in \Phi_{x}(z_{t},\lambda_{0})$, the previous inequality yields
	\[
	0 = t \langle z^{\ast}_{t}, y - z_{t} \rangle + (1-t) \langle z^{\ast}_{t}, x - z_{t} \rangle \leq t \langle z^{\ast}_{t}, y - z_{t} \rangle
	\]
	so that $\inf_{z^{\ast}_{t} \in \Phi_{x}(z_{t})} \langle z^{\ast}_{t}, y - z_{t} \rangle \geq 0$ for $t \in ]0,1[$. The upper sign-continuity of $\Phi_{x}(\cdot,\lambda_{0})$ then implies
	\[
	\sup_{x^{\ast} \in \Phi_{x}(x,\lambda_{0})} \langle x^{\ast}, y - x \rangle \geq 0,
	\]
	where the $w^{\ast}$-compactness of $\Phi_{x}(x,\lambda_{0})$ guarantees that the supremum is attained at some $x^{\ast}_{0} \in \Phi_{x}(x,\lambda_{0}) \subset T(x,\lambda_{0})\setminus\{0\}$. This shows the closedness of $LSVI^{w,int}_{pert}$.
	
	Now, assume that $\Phi_{x}(x,\lambda_{0})$ is convex.  As $x \in LSVI^{w,int}_{pert}(\lambda_{0},\mu_{0})$, we have
	\[
	\inf_{y \in V_{x}\cap K(\mu_{0})} \; \sup_{x^{\ast} \in \Phi_{x}(x,\lambda_{0})} \langle x^{\ast}, y - x \rangle \geq 0.
	\]
	The Sion minimax theorem then guarantees the existence of an element $x^{\ast} \in LSVI^*_{pert}(\lambda_{0},\mu_{0})$. Thus the closedness of $LSVI^*_{pert}$ is obtained.\qed
\end{proof} 

\begin{remark}\label{rmk:beta_partial}
	Looking into the above proof, we may deduce the following partial results which will be useful for us in the sequel.
	
	Suppose that $(\lambda_{n},\mu_{n})_n$ is a sequence in $\Lambda \times U$ converging to $(\lambda_{0},\mu_{0}) \in \Lambda \times U$ such that \ref{assmp:loc_upper_sign} of Theorem \ref{thm:beta} and that the following conditions hold:
	\begin{enumerate}[label = \upshape \roman*')]
		\item $LSVI^{w,int}_{pert}(\lambda_{n},\mu_{n}) \neq \emptyset$, for all $n\in\N$;
		\item $K(\mu_{n})$ is locally convex with nonempty interior, for all $n\in\N$;
		\item the couple of set-valued maps $(T\setminus\{0\},K)$ is int-dually lower semicontinuous on $X^2\times\Lambda\times U$;
		\item $(K(\mu_{n}))_n$ is Mosco convergent to $K(\mu_{0})$.
	\end{enumerate}
	If $x_{n} \in LSVI^{w,int}_{pert}(\lambda_n,\mu_n)$ for all $n\in\N$ and $(x_{n})_n$ converges to $x$, then $x \in LSVI^{w,int}_{pert}(\lambda_0,\mu_0)$. Furthermore, if $\Phi_{x_{0}} (x_{0},\lambda_{0})$ is convex, then $x_{0} \in LSVI^*_{pert} (\lambda_{0},\mu_{0})$.
\end{remark}

\section{Applications to optimization}\label{sec:optim}

In this last section, our aim is to explore the applications to optimization of the above developed local analysis for variational and quasivariational inequalities. Our final target being, in Subsection \ref{subsec:SLLMF}, the new concept of Single-Leader-Local-Multi-Follower game, that is hierarchical games in which the equilibrium between the followers is considered in a local sense, we will first consider, respectively in subsection \ref{subsec:optim} and \ref{subsec:quasi-optim}, the case of perturbed optimization and quasi-optimization problems.

\subsection{Perturbed constrained optimization problems}\label{subsec:optim}

Let us now consider the solution map $LOpt$ of a perturbed optimization problem defined, for any $(\lambda,\mu)\in \Lambda\times U$, by $LOpt(\lambda,\mu)$ being the set of local solutions (in the classical sense) of 
\[\begin{array}{rl}
	\min_x & f(x,\lambda)\\[3pt]
	\mbox{s.t.} & x\in K(\mu).
	\end{array}
\]
If $f$ is quasiconvex in its first argument, then we abuse the notation of $F_{f}$ (as defined in \eqref{eqn:Ff_mapping}) to be understood as
\[
	F_{f} (x,\lambda) := F_{f(\cdot,\lambda)}(x)
\]
for $(x,\lambda) \in X \times \Lambda$. Thus, for example, saying that $(F_{f},K)$ is int-dually lower semicontinuous is to be interpreted in the above sense.

\begin{theorem}\label{thm:beta_VI_f}
	Suppose that $X, \Lambda, U$ are all finite-dimensional, $f : X \times \Lambda \to \R$ is lower semicontinuous in both variables, and is continuous and semistrictly quasiconvex on the first argument. Assume that the following conditions are fulfilled:
	\begin{enumerate}[label = \upshape \roman*)]
		\item\label{assmp:level_bdd} $f(x,\lambda)$ is level-bounded in $x$ locally uniformly in $\lambda$.
		\item\label{assmp:cont-cont} at each $\lambda \in \Lambda$, there is a point $x_{\lambda} \in \argmin_{X} f(\cdot,\lambda)$ such that $f(x_{\lambda},\cdot)$ is continuous on $\Lambda$.
		\item for all $\mu \in U$, $K(\mu)$ is locally convex with nonempty interior;
		\item the couple {$(F_f\setminus\{0\},K)$} is int-dually lower semicontinuous $X^2\times\Lambda\times U$;
		\item for any sequence $(\mu_{n})_n$ converging to $\mu$, the sequence $(K(\mu_{n}))_{n}$ is Mosco convergent to $K(\mu)$.
	\end{enumerate}
	Then, the map $(\lambda,\mu) \mapsto LOpt(\lambda,\mu)$ is closed.
\end{theorem}
The following proposition will be useful in the proof of the above theorem.
\begin{proposition}\label{prop:iff_quasiconvex_optim}
Let $f:X\to\R$ be sub-boundarily constant and semistrictly quasiconvex on $X$ and $C$ be a nonempty locally convex subset of $X$. Then $\bar{x}$ is a local minimum of $f$ over $C$ if and only if it solves $LSVI(F_{f},C)$. Moreover, $LSVI(F_{f},C)$ is a union of two disjoint sets $\argmin_{X} f$ and $ LSVI^{*}(F_f,C) := LSVI(F_{f}\setminus\{0\},C)$ {, that is 
\[
	\textstyle\argmin_{C}^{loc} f = LSVI(F_{f},C) = \left[\argmin_{X} f \cap C \right] \cup LSVI^{*}(F_f,C).
\]
where, as usual, {$\argmin^{loc}_{C} f$} stands for the set of local minumum of $f$ on $C$.
}
\end{proposition}
\begin{proof}
	First note that $\argmin_{X} f \cap C$ and $ LSVI^{*}(F_f,C) := LSVI(F_{f}\setminus\{0\},C)$ are disjoint by Proposition \ref{prop:Ff}. The union of them being $LSVI(F_{f},C)$ is obvious. Let us now prove the sufficient condition only for the case that $\bar{x} \in LSVI^{*}(F_{f},C)$ which means that there exists a neighbourhood $U_{\bar{x}}$ of $\bar{x}$ such that $\bar{x}\in SVI(F_{f}\setminus\{0\},C\cap U_{\bar{x}})$. Without loss of generality $U_{\bar{x}}$ can be assumed to be convex. Thus, by \cite[Proposition 2.9]{A_CoaVan_Salas21}, $\bar{x}$ is a local minimum of $f$ over $C$. Observing that, since $f$ is {sub-boundarily} constant on $X$, the sublevel set $S_{f(\bar{x})}$ has nonempty interior, then the necessary conditions can be {deduced} in the same lines as \cite[Corollary 5.2]{A_chapter14}.     \qed
\end{proof}

\DA{The above} formula highlights the necessary and sufficient local optimality condition in term of variational inequalities, the same relation does not extend to the global setting. Let us illustrate this fact in the following simple example.
\begin{example}
	Consider $X = \R$, $C = [-2,-1] \cup [1,2]$ and define $f : X \to \R$ by $f(x) = \abs{x}$. Observe first that the global and local minimizers over $C$ are identical with $\argmin_{C} f = \argmin_{C}^{loc}= \{-1,1\}$. On the other hand, we have
	\[
		F_{f} (x) = \begin{cases}
			\{-1\}	&\text{if $x < 0$}\\
			[-1,1]	&\text{if $x = 0$}\\
			\{1\}	&\text{if $x > 0$}
		\end{cases}
	\]
	which yields $LSVI(F_{f},C) = LSVI^{*}(F_{f},C) = \argmin_{C}^{loc} = \{-1,1\}$. However, no global solution to $SVI(F_{f},C)$ exists. We conclude in this example that the global solutions of $\argmin_{C}f$ cannot be expressed with that of $SVI(F_{f},C)$. 

	\blue{The above conclusion is actually to be expected, since $SVI(F_{f},C)$ generally fails to act as a necessary optimality condition when $C$ is not convex. On the other hand, the conditions of $C$ and $f$ here allows us to write a necessary optimality condition due to Aussel and Ye \cite[Theorem 4.1]{zbMATH05132179} for local minimizers in $\argmin_{C}^{loc} f$ as 
	\[
		\textstyle\bar{x} \in \argmin_{C}^{loc} f \implies 0 \in N_{f}(\bar{x}) + N_{L}(C;\bar{x}),
	\]
	where $N_{L}(C;\bar{x})$ denotes the limiting normal cone of $C$ at $\bar{x}$. \DA{Note} that $SVI(F_{f},C)$ also generally fails as a sufficient optimality condition without the convexity of $C$ (see \cite[Proposition 3.2]{zbMATH05132179}).
	} \qedex
\end{example}

\begin{proof}{\sl of Theorem \ref{thm:beta_VI_f}.}~
Taking into account the semistrict quasiconvexity of $f$, let us first define the set-valued map $LSVI^*_{pert}$: for any $(\lambda,\mu)\in \Lambda\times U$, let $LSVI^*_{pert}(\lambda,\mu)$ be the set of $\bar{x}\in K(\mu)$ for which there exist a neighbourhood $U_{\bar{x}}$ of $\bar{x}$ and $\bar{x}^{\ast}\in F_f(\bar{x},\lambda)$ with $\bar{x}^{\ast}\neq 0$ such that 
	\[\langle \bar{x}^*, y-\bar{x}\rangle \ge 0,\quad\forall\,y\in K(\mu)\cap U_{\bar{x}}.
	\]
In other words for any $(\lambda,\mu)\in \Lambda\times U$, 
\[LSVI^*_{pert}(\lambda,\mu)=LSVI(F_{f(\cdot,\lambda)}\setminus\{0\},K(\mu)).
\]
Thus combining the equality with the fact that any continuous function is {sub-boundarily} constant, one can deduce from Proposition \ref{prop:iff_quasiconvex_optim} that 
\[LOpt(\lambda,\mu) = LSVI^*_{pert}(\lambda,\mu) \cup \left[\argmin_{X} f(\cdot,\lambda) \cap K(\mu) \right],\quad\forall (\lambda,\mu) \in \Lambda \times U.
\]
Now, let $(\lambda_{n},\mu_{n})_n$ and $(x_{n})_n$ be two sequences in $\Lambda \times U$ and $X$ converging respectively to $(\lambda_0,\mu_0)$ and $x$ with
\[x_{n} \in {LOpt(\lambda_{n},\mu_{n})},\quad\mbox{ for all }n \in \N.
\]
If $(x_{n})_n$ contains a subsequence $(x_{n_{k}})_k$ such that $x_{n_{k}} \in LSVI^*_{pert}(\lambda_{n_{k}},\mu_{n_{k}})$ for all $k \in \N$, then Remark \ref{rmk:beta_partial}{, combined with Proposition \ref{prop:iff_quasiconvex_optim},} shows that $x \in LSVI^*_{pert}(\lambda_0,\mu_0) \subset LOpt(\lambda_0,\mu_0)$. 
Otherwise, $x_{n} \in \left[\argmin_{X} f(\cdot,\lambda_{n}) \cap K(\mu_{n}) \right]$ for all sufficiently large $n \in \N$. By \cite[Theorem 1.17 (c) ]{zbMATH01110799}, $\inf_{X} f(\cdot,\lambda_{n})$ converges to $\inf_{X}f(\cdot,\lambda_{0})$. Then the proof is complete since, according to \cite[Theorem 7.41 (a)]{zbMATH01110799} and the Mosco convergence of $K(\mu_{n})$ to $K(\mu_{0})$, we get $x \in \left[\argmin_{X} f(\cdot,\lambda_{0})\cap K(\mu_{0})\right] \subset LOpt(\lambda_0,\mu_0)$. \qed
\end{proof}
The closedness of the solution map, in the global sense, of perturbed optimization problems (and variational inequalities) has been extensively study in the literature (see e.g. \cite{zbMATH01502618}). On the contrary, the literature concerning the solution map in the local sense is rather scarce. See e.g. \cite{zbMATH01417166,zbMATH06574936} for the few existing {resources} in this direction. In this section, Theorem \ref{thm:beta_VI_f} provides, in the case of semistrictly quasiconvex objective functions, alternative sufficient conditions for the closedness of the solution map in the local sense.

\subsection{{Perturbed quasi-optimization problems}}\label{subsec:quasi-optim}
{A quasi-optimization problem is an optimization problem in which the constraint set depends on the considered point. \blue{This terminology has been introduced in \cite{zbMATH01568702} and existence results for such problems has been then obtained\blue{, in the global sense,} in \cite{A_Cotrina13} by the use of the normal operator $N^a$ to adjusted sublevel sets. We quote also the recent work of in \cite{zbMATH06780576} based on the direct approach leading to the so-called $\lambda$-eigenvalue solutions, see \cite[Th. 9.9]{zbMATH06780576} as well as \cite[Remark 9.10]{zbMATH06780576} for a deep discussion and comparison with the similar results of \cite{A_Cotrina13}.} 

Our aim in this subsection is to make an analysis of local solutions of quasi-optimization problem: reformulation, existence and closedness will be investigated respectively in the forthcoming Proposition \ref{prop:LQOpt_as_union}, Theorem \ref{exist_quasi_opt} and Theorem \ref{thm:gamma}.
}

{But let us first precise the concept of local solution in case of (perturbed or not) quasi-optimization problems.
\begin{definition}\label{def:QOpt}$~$\\[-10pt]
\begin{itemize}
	\item Consider a function $f:X \to \R$ and a set-valued map $K:X\multimap X$. A point $\xb\in X$ is called a local solution of the  \emph{quasi-optimization problem} defined by $f$ and $K$ if $\xb\in K(\xb)$ and there exists a neighbourhood $U_{\xb}$ of $\xb$ such that 
\begin{equation}\label{pb_quasi_opt_local}
	f(\xb)\le f(x), \quad\forall\,x\in K(\xb)\cap U_{\xb}.
\end{equation}
The set of such local solutions will be denoted by $LQOpt(f,K)$.
	\item Consider a function $f:X\times \Lambda \to \R$ and a set-valued map $K:X\times U\multimap X$. For any $(\lambda,\mu)\in \Lambda\times U$, a point $\xb\in X$ is called a local solution of the  \emph{perturbed quasi-optimization problem} defined by $f(\cdot,\lambda)$ and $K(\cdot,\mu)$ if $\xb\in K(\xb,\mu)$ and there exists a neighbourhood $U_{\xb,\mu}$ of $\xb$ such that 
\begin{equation}\label{pb_quasi_opt_local}
	f(\xb,\lambda)\le f(x,\lambda), \quad\forall\,x\in K(\xb,\mu)\cap U_{\xb,\mu}.
\end{equation}
The set of such local solutions will be denoted by $LQOpt(f(\cdot,\lambda),K(\cdot,\mu))$.
\end{itemize}
\end{definition}
}

The following result reveals {relationships} between solution sets of the problems LOpt, LQOpt and LSVI {or, more precisely provides us with two reformulations of the quasi-optimization problem, understood in the local sense.} 

\begin{proposition}\label{prop:LQOpt_as_union}
	Let us consider a perturbed function $f:X\times \Lambda \to \R$ and a perturbed set-valued map $K:X\times U\multimap X$. Let us assume that, for any $\mu\in U$, the partial set-valued map $K(\cdot,\mu)$ is locally reproducible on $FP(K(\cdot,\mu))$ (with $U_{z,\mu}$ being a neighbourhood of reproducibility of $K(\cdot,\mu)$ at $z$). Then{, for any $(\lambda,\mu)\in\Lambda\times U$,}
	\[
	LQOpt(f(\cdot,\lambda),K(\cdot,\mu)) = \bigcup_{z\in FP(K(\cdot,\mu))} LOpt(f(\cdot,\lambda),{K(z,\mu)\cap U_{z,\mu}}).
	\]
	Moreover, if $f$ is continuous and semistrictly quasiconvex in the first argument, then we also have 
	\[
	LQOpt(f(\cdot,\lambda),K(\cdot,\mu)) = \bigcup_{z \in FP(K(\cdot,\mu))} LSVI(F_{f(\cdot,\lambda)},K(z,\mu)\cap U_{z,\mu}).
	\]
\end{proposition}
\begin{proof}
	Let us first assume that $z$ is a fixed point of the partial set-valued map $K(\cdot,\mu)$ and $x\in LOpt(f(\cdot,\lambda),{K(z,\mu)\cap U_{z,\mu}})$. Thus there exists a neighbourhood $U_x$ of $x$ such that $x\in {K(z,\mu)\cap U_{z,\mu}\cap U_x}$ and 
	\begin{equation}
	f(x)\le f(u),\qquad \forall\, u\in K(z,\mu)\cap U_{z,\mu}\cap U_x.
	\label{eq:local_formula}
	\end{equation}
	Since $K(\cdot,\mu)$ is locally reproducible at $z$ and $x\in K(z,\mu)\cap U_{z,\mu}$, one immediately has {$K(x,\mu)\cap U_{z,\mu}=K(z,\mu)\cap U_{z,\mu}$}. Therefore $x\in K(x,\mu)\cap U_{z,\mu}\cap U_x$ and combining with \eqref{eq:local_formula}, $x\in LQOpt(f(\cdot,\lambda),K(\cdot,\mu))$.
	
	The inverse inclusion is trivial since any element $x$ of $LQOpt(f(\cdot,\lambda),K(\cdot,\mu))$ is a fixed point of the partial map $K(\cdot,\mu)$ and also element of the local minimizer set $LOpt(f(\cdot,\lambda),{K(x,\mu)\cap U_{x,\mu}})$. 
	
{The second formula can be deduced from Proposition \ref{prop:iff_quasiconvex_optim}.}
\qed
\end{proof}
Combining the Proposition \ref{prop:LQOpt_as_union} above with a classical existence result for quasimonotone Stampacchia variational inequalities, one can obtain sufficient conditions for the existence of local solutions of the perturbed quasi-optimization problem \eqref{pb_quasi_opt_local}.

\begin{theorem}\label{exist_quasi_opt}
	Let us assume that{, for any $\lambda\in \Lambda$,} the function {$f(\cdot,\lambda)$} is quasiconvex and sub-boundarily constant, and {that, for any $\mu\in U$, the partial map $K(\cdot,\mu)$} is locally reproducible at some fixed point $\bar{z} \in {FP(K(\cdot,\mu))}$.
	If {$K(\cdot,\mu)$} has nonempty locally convex compact values on $X$, then {the perturbed quasi-optimization problem defined by $f(\cdot,\lambda)$ and $K(\cdot,\mu)$ admits a local solution, that is $LQOpt(f(\cdot,\lambda),K(\cdot,\mu))\neq\emptyset$.}
\end{theorem}
\begin{proof}
	Let $z \in {FP(K(\cdot,\mu))}$ and {$U_{z,\mu}$ be} a neighbourhood of reproducibility at $z$ such that {$K(z,\mu)\cap U_{z,\mu}$} is convex and compact. {If $\argmin_{X} f(\cdot,\lambda)\cap [K(z,\mu)\cap U_{z,\mu}]\neq\emptyset$, there is nothing to prove. Otherwise, according to Proposition \ref{prop:Ff}, $F_f$ is locally upper sign-continuous (in the first variable) and $0\not\in F_{f(\cdot,\lambda)}$ on $K(z,\mu)\cap U_{z,\mu}$. Then \cite[Theorem 2.1]{zbMATH02125839} ensures the nonemptiness of {the solution set} $LSVI(F_{f},{K(z,\mu)\cap U_{z,\mu}})=LSVI(F_{f}\setminus\{0\},{K(z,\mu)\cap U_{z,\mu}})$, and so is $LOpt({f(\cdot,\lambda),K(z,\mu)\cap U_{z,\mu}})$ since, by \cite[Proposition 2.9]{A_CoaVan_Salas21}, the set $LSVI(F_{f}\setminus\{0\},{K(z,\mu)\cap U_{z,\mu}})$ is included in $LOpt({f(\cdot,\lambda),K(z,\mu)\cap U_{z,\mu}})$.
}	
\qed
\end{proof}

\blue{Finally let us present the qualitative stability for the perturbed quasi-optimization problem \eqref{pb_quasi_opt_local} (closedness property on the local solution map) which is a natural extension of the closedness results proved by Aussel and Sagratella \cite{A_Sagratella17}, Ait Mansour and Aussel \cite{zbMATH05351708} and Aussel and Cotrina \cite{A_Cotrina13}.}

\begin{theorem}\label{thm:gamma}
	Suppose that $X, \Lambda, U$ are all finite-dimensional, $f : X \times \Lambda \to \R$ is continuous and semistrictly quasiconvex in the first argument {and lower semicontinuous in the second argument}. Let $K : X \times U \multimap X$ be a closed {set-valued} map whose values have nonempty interior and are locally convex. Suppose that \ref{assmp:level_bdd} and \ref{assmp:cont-cont} of Theorem \ref{thm:beta_VI_f} {hold, that} $FP(K(\cdot,U))$ is compact and {that,} for each $(\lambda,\mu) \in \Lambda \times U$ and $z \in FP(K(\cdot,U))$, the following conditions hold:
	\begin{enumerate}[label=\upshape\roman*)]
			\item $K(\cdot,\mu)$ is locally reproducible on $FP(K(\cdot,\mu))$ (with $U_{z,\mu}$ being the neighbourhood of reproducibility at $z$);
			\item {the couple of maps $(F_f\setminus\{0\},K)$ is int-dually lower-semicontinuous on $X^2\times \Lambda\times U$};
			\item for any {sequence $(\mu_{n})_n$ converging to} $\mu$, {the sequence $(K(z,\mu_{n}))_n$} is Mosco convergent to $K(z,\mu)$.
	\end{enumerate}
	Then the map $(\lambda,\mu) \mapsto LQOpt(f(\cdot,\lambda),K(\cdot,\mu))$ is closed.
\end{theorem}

\blue{Before continuing to the required technical lemmas and proof of this theorem, we would like to make the following remarks with respect to the global setting.
\begin{remark}\label{rmk:LQOpt_closedness}~\\[-1,2em]
	\begin{enumerate}[label=\upshape\alph*)]
		\item Theorem \ref{thm:gamma} covers the global case when the reproducibility condition and solutions are relaxed to global ones.
		\item\label{rmk:LQOpt_others} In the global treatment of stability in \cite{Ait-Mansour_A06}, the upper-semicontinuity of the solution map is obtained from quantitative stability (see \cite[Corollary 3.3]{Ait-Mansour_A06}). For global parametric quasi-optimization problems, a similar conclusion of qualitative stability is obtained, again, from the quantitative approach in \cite[Proposition 3.5]{zbMATH07538430}.
	\end{enumerate}
\end{remark}
}

In order to deduce the closedness for the local solution map of the perturbed quasi-optimization problem {\eqref{pb_quasi_opt_local}}, we need the following technical lemmas.
	
\begin{lemma}\label{lem:USC-union}
	Suppose that $C \subset X$ is closed and that the set-valued maps $G : C \multimap C$ and $T : C \times X \multimap C$ are upper semicontinuous respectively on $C$ and $C\times X$, both having compact values. Then the map $F : C \multimap C$ given, for any $x\in C$, by
	\[
		F(x) := \bigcup_{u \in G(x)} T(x,u)
	\]
	is upper semicontinuous. In particular, $F$ is closed.
\end{lemma}

\begin{proof}
	Fix $x \in C$ and pick an arbitrary $\varepsilon > 0$. Since, for any $u\in G(x)$, $B(F(x),\varepsilon)$ is a neighbourhood of $T(x,u)$, the upper semicontinuity of $T$  implies that, for each $u \in G(x)$, there exist $\delta_{u},\eta_{u} > 0$ such that $T(x',u') \subset B(F(x),\varepsilon)$ for all $(x',u') \in \gr G \cap (B(x,\delta_{u}) \times B(u,\eta_{u}))$. Since $G(x)$ is compact, it can be covered by a finite family of open balls $\{B(u_{i},\eta_{u_{i}})\}_{i=1}^{p}$ for some $u_{1},\cdots,u_{p} \in G(x)$. By the upper semicontinuity of $G$, there exists $\delta > 0$ such that $G(\hat{x}) \subset \bigcup_{i=1}^{p} B(u_{i},\eta_{u_{i}})$ whenever $\hat{x} \in B(x,\delta) \cap C$. Put $\delta_{0} := \min\{\delta,\delta_{u_{1}},\cdots,\delta_{u_{p}}\}$, we have $F(x') \subset B(F(x),\varepsilon)$ for all $x' \in B(x,\delta_{0})$.   \qed
\end{proof}
	
\begin{lemma}\label{lem:USC-FP}
	Suppose that $C \subset X$ is closed and $T : C \times X \multimap C$ is closed, then the fixed point map $G(z):= FP(T(\cdot,z))$ (defined for each $z \in C$) is closed.
\end{lemma}

\begin{proof}
		Let $(x_{n},y_{n})_n$ be a sequence in $\gr G$ which converges to a point $(x,y) \in X \times C$. This means for any $n \in \N$, we have $y_{n} \in T(y_{n},x_{n})$ and hence $(y_{n},x_{n},y_{n}) \in \gr T$. Since {$T$} is closed, $y \in T(y,x)$ so that $(x,y) \in \gr G$.   \qed
\end{proof}
	
\begin{proof}{\em of Theorem \ref{thm:gamma}} 
Let $z \in FP(K(\cdot,U))$. {Combining propositions \ref{prop:prop_normal} and \ref{prop:Ff} with the hypothesis that the sets $K(z,\mu)$ are locally convex, one can deduce from Theorem \ref{thm:existence_LSVI} that, for any $(\lambda,\mu)$, $LSVI(F_{f(\cdot,\lambda)},K(z,\mu)\cap U_{z,\mu})$ is nonempty. On the other hand, since the couple {$(F_f\setminus\{0\},K(z,\mu))$} is int-dually lower semicontinuity, so is the couple {$(F_f\setminus\{0\},K(z,\mu)\cap U_{z,\mu})$} and thus, according to Theorem \ref{thm:beta_VI_f}, the map ${(\lambda,\mu)\mapsto} LOpt(f(\cdot,\lambda),{K(z,\mu)}\cap U_{z,\mu})$ is closed. By Proposition \ref{prop:iff_quasiconvex_optim}, the map ${(\lambda,\mu)\mapsto} LSVI(F_{f(\cdot,\lambda)}, K(\cdot,\mu)\cap U_{z,u})$ is also closed and since the latter is true for any $z \in FP(K(\cdot,U))$, Proposition \ref{prop:LQOpt_as_union}, together with the local reproducibility of the map $K(\cdot,\mu)$ on $FP(K(\cdot,\mu))$,  allows to conclude to the closedness of the map ${(\lambda,\mu)\mapsto} LQOpt(f(\cdot,\lambda),K(\cdot,\mu))$.
}	\qed
\end{proof}

\subsection{{Single-Leader-Local-Multi-Follower games}}\label{subsec:SLLMF}

{A Single-Leader-Multi-Follower game (SLMF) corresponds to a model involving $M+1$ players in a hierarchical/non-cooperative interaction. More precisely $M$ of the players (called {\em followers}), each deciding a variable $y_i$, try to reach a non-cooperative generalized Nash equilibrium, while their optimization problem is parametrized by the decision $y_{-i}=(y_1,\dots,y_{i-1},y_{y+1},\dots, y_M)$ and by the decision $x$ of player $M+1$, called the {\em leader}. In the {\em optimistic formulation} of a (SLMF) game, the leader minimizes his objective function with respect to both variables $x$ and $y=(y_1,\dots,y_M)$. Single-Leader-Follower games has been extensively studied and allow to model numerous applications in energy management, transport, economics, see \cite{A_Svensson_chapter21} and references therein.
}

{
Our aim in this last section is to define and analyse a slightly modified version of (SLMF) problem in which the generalized Nash equilibrium $y$ obtained by the followers is actually considered in the local sense, thus leading to a {\em Single-Leader-Local-Multi-Follower} (SLLMF) game.
}

{In the context of a non-cooperative interactions parametrized by the leader variable $x$, let us first recall the concept of local Nash equilibrium. So for any $i=1,\dots,M$, let $y_i\in \R^{m_i}$ be the decision variable of player $i$, $\theta_i(y_i,y_{-i},x)$ be his objective function and $K_i(y_{-i},x)$ be his constraint set. Then a point $y=(y_1,\dots,y_M)$ is a local generalized Nash equilibrium if, for any $i=1,M$, there exists a neighbourhood $U_i$ of $y_i$ such that $y_i$ solves the following parametrized problem 
\begin{equation}\label{def_local_Nash}
	\begin{array}{rl}
		\min_{z_i} & \theta_i(z_i,y_{-i},x)\\[5pt]
		\mbox{s.t.} & z_i\in K_i(y_{-i},x)\cap U_{i}.
	\end{array}
\end{equation}
We denote by $LGNEP(\theta_i(\cdot,x),K_i(\cdot,x))_{i=1}^M$ the set of all such local generalized Nash equilibria.
} \blue{This concept of local generalized Nash equilibrium was also considered by Rockafellar in \DA{\cite{zbMATH06836808,Terry_preprint}}}. 

{Now given two nonempty sets $C_{1} \subset \R^{n}$, $C_{2} \subset \R^{m}$, with $m=\sum_i m_i$ and the objective function $F : \R^{n} \times \R^{m} \to \R$ of the leader, the (optimistic) Single-Leader-Local-Multi-Follower games consists in finding a couple $(x,y)\in \R^{n} \times \R^{m}$ solution of 
\[\begin{array}{rrl}
	(SLLMF)\qquad\qquad&\min_{x,y} & F(x,y)\\
	&\text{s.t.}& \left\{
	\begin{array}{l}
	x\in C_1\\
	y\in C_2\\
	y\in LGNEP(\theta_{i}(\cdot,x),K_{i}(\cdot,x))_{i=1}^{M}.
	\end{array}
	\right.
	\end{array}
\]
One of the main motivation to consider Single-Leader-Multi-Follower model in which the Nash equilibrium is understood in the local sense is that it is well-known that computing a generalized Nash equilibrium is not an easy task and that this computation is often replaced by the resolution of the associated {\em concatenated Karush-Kuhn-Tucker system}, thus leading, if convexity assumptions are not satisfied, to local generalized Nash equilibria. {Interested readers may consult, e.g. \cite{Facchinei_Pang_book1,zbMATH05995776}, for references.}
}

{Following a classical approach (see e.g. \cite{Facchinei_Pang_book1}), the use of the so-called {\em Nikaido-Yosida} function allows to describe the set of equilibria of a generalized Nash equilibrium game as the zero (and minimums) of an associated {\em gap function}. The same can be done for local equilibrium and the aim of the following theorem is to adapt it to the context of interactions between the followers of the (SLLMF) game. Its proof is omitted and can be developed in line with \cite{Facchinei_Kanzow07}, but stated in the local version.
}

\begin{theorem}\label{thm:NI}
{Let $x$ be an element of $\R^n$ and suppose that $\Psi$ is the Nikaido-Isoda function associated to the objective function $\theta_i(\cdot,x)$, $i=1,m$ of the followers, that is
\begin{equation}\label{eqn:NI_function}
	\Psi(y,z,x) := \sum_{i = 1}^{M} [ \theta_{i}(y_{i},y_{-i},x) - \theta_{i}(z_{i},y_{-i},x) ],\quad\forall\,(y,z)\in \R^m\times \R^m.
\end{equation}
Define, for each $y \in FP(K(\cdot,x))$ and each neighbourhood of reproducibility $U_{y,x}$ of $K(\cdot,x)$ at $y$, the gap function 
\[{V^{U_{y,x}}}(y,x) := \sup_{z \in K(y,x) \cap U_{y,x}} \Psi(y,z,x).
\]
Then, if the product map $K(y,x) := \prod_{i=1}^{M}K_{i}(y_{-i},x)$ is locally reproducible at each $y \in FP(K(\cdot,x))$, the following properties hold:
\begin{enumerate}[label=\upshape\alph*)]
		\item \label{conclusion:NI-a} ${V^{U_{y,x}}}(y,x) \geq 0$ for all $y \in FP(K(\cdot,x))$ and all choices of $U_{y,x}$;
		\item \label{conclusion:NI-b} If $\bar{y} \in FP(K(\bar{y},x))$, then ${V^{U_{\bar{y},x}}}(y,x))$ is defined for all $y \in FP(K(\bar{y},x)) \cap U_{\bar{y},x}$;
		\item \label{conclusion:NI-c} $\bar{y} \in FP(K(\bar{y},x))$ and ${V^{U_{\bar{y},x}}}(\bar{y},x) = 0$ if and only if $\bar{y}$ is an element of $LGNEP(\theta_{i}(\cdot,x),K_{i}(\cdot,x))_{i=1}^{M}$ with the same neighbourhood;
		\item \label{conclusion:NI-d}   If $LGNEP(\theta_{i}(\cdot,x),K_{i}(\cdot,x))_{i=1}^{M}$ is nonempty and $U_{y,x}$ is fixed at each solution $y \in LGNEP(\theta_{i}(\cdot,x),K_{i}(\cdot,x))_{i=1}^{M}$, then it holds 
		\[
			LGNEP(\theta_{i}(\cdot,x),K_{i}(\cdot,x))_{i=1}^{M} = LQOpt(V(\cdot,x),K(\cdot,x)).
		\]
		where $V(y,x)$ stands for $V(y,x) := {V^{U_{y,x}}}(y,x)$.
\end{enumerate}
}
\end{theorem}
Our aim is here to prove the theorem below proposing a set of assumptions under which the Single-Leader-Local-Multi-Follower game (SLLMF) admits at least a solution. 

\begin{theorem}\label{thm:exist_SLMFG}
	Consider a $(SLLMF)$ with $F$ being {lower semicontinuous}, $\theta_{i}$ being continuous for all $i = 1,\cdots,M$. Suppose that $C_{1}$ and $C_{2}$ are closed sets, that $y \in C_{2} \mapsto K(y,x)$ is graph-compact and locally reproducible over its fixed point set for each $x \in C_{1}$, that  $K$ is a lower semicontinuous set-valued map in the second argument and that $FP(K(\cdot,C_{1}))$ is closed. If $V$ is semistrictly quasiconvex and the couple $(F_{V},K)$ is int-dually lower semicontinuous on $C_{2}^{2} \times C_{1} \times C_{1}$, where $V$ is defined as in Theorem \ref{thm:NI}, then the (SLLMF) admits at least a solution.
\end{theorem}

Some of the hypotheses of this existence result, like the int-dual lower semicontinuity of $(F_{V},K)$, can appear to be quite restrictive. As shown in Theorem \ref{thm:exist_SLMFG_linear}, this confinement condition can be omitted when, for example, all $\theta_{i}$'s are linear and the constraint maps $K_{i}$'s are linear translations. It is also important to notice that, as observed in \cite{A_Svensson_chapter21}, there are very few existence result for {Multi-Leader-Multi-Follower} problems (see e.g. historical comments in \cite[Section 3.4.1]{A_Svensson_chapter21}). In the case of {Single-Leader-Multi-Follower} games a general result has been proved in \cite{A_Svensson_chapter21} (see Theorem 3.3.4). The proof of this existence result and the one of Theorem \ref{thm:exist_SLMFG} and the forthcoming Theorem \ref{thm:BOpt_existence} are very similar but the assumptions slightly differ: indeed while here we have a stronger hypothesis of int-dual lower semicontinuity of $(F_{V},K)$, the quite restrictive assumptions on the lower semi-continuity (with respect to both variables $x$ and $y_{-i})$ made in \cite[Theorem 3.3.4]{A_Svensson_chapter21} on the constraint map $K_i$ is replaced here by the lower semi-continuity of the constraint map $K_i$ with respect to the second argument (the leader variable $x$) coupled with its local reproducibility. Both existence results \cite[Theorem 3.3.4]{A_Svensson_chapter21} and the above stated Theorem \ref{thm:exist_SLMFG} are complementary.

The proof of Theorem \ref{thm:exist_SLMFG} relies on the reformulation, thanks to the Nikaido-Isoda gap function, of the (SLLMF) problem as a \emph{bilevel local quasi-optimization problem} (briefly, BLQOpt). Let us start by precising the meaning of (BLQOpt). Let $C_{1} \subset \R^{n}$ and $C_{2} \subset \R^{m}$ be nonempty, $F : \R^{n} \times \R^{m} \to \R$, $f : \R^{m} \times \R^{n} \to \R$ and $K : C_{2} \times C_{1} \multimap C_{2}$. The associated (BLQOpt) is given as
\[
	\begin{array}{rrl}
	(BLQOpt)\qquad\qquad&\min_{x,y} & F(x,y)\\
	&\text{s.t.}& \left\{
	\begin{array}{l}
	x\in C_1\\
	y\in C_2\\
	y\in LQOpt(f(\cdot,x), K(\cdot,x)).
	\end{array}
	\right.
	\end{array}
\]

Thus, as an intermediate step to the proof of Theorem \ref{thm:exist_SLMFG}, let us give sufficient conditions for the existence of solutions of a (BLQOpt).

\begin{theorem}\label{thm:BOpt_existence}
	Suppose that $C_{1} \subset \R^{n}$ and $C_{2} \subset \R^{m}$ are nonempty and closed, and $F : \R^{n} \times \R^{m} \to \R$ is {lower semicontinuous}. Also suppose that $f : \R^{m} \times \R^{n} \to \R$ is a function in which for each $x \in C_{1}$, $f(\cdot,x)$ {is} level-bounded on $C_{2}$, continuous and semistrictly quasiconvex. Let $K : C_{2} \times C_{1} \multimap C_{2}$ be a graph-compact map which is lower semicontinuous in the second argument and its values have nonempty interior and are locally convex. Suppose that $FP(K(\cdot,C_{1}))$ is closed and for each $x \in C_{1}$ and $z \in FP(K(\cdot,x))$, the following conditions hold:
	\begin{enumerate}[label=\upshape\roman*)]
		\item $K(\cdot,x)$ is locally reproducible on its fixed point set $FP(K(\cdot,x))$;
		\item\label{hyp2_normalized} the couple of maps $(F_f\setminus\{0\},K)$ is int-dually lower semicontinuous on $C_{2}^{2} \times C_{1} \times C_{1}$.
	\end{enumerate}
	If the graph of $x \mapsto LQOpt(f(\cdot,x),K(\cdot,x))$ is nonempty, then $(BLQOpt)$ {admits at least} a solution.
\end{theorem}
\begin{proof}
	The map $x \mapsto LQOpt(f(\cdot,x),K(\cdot,x))$ is graph-closed by Theorem \ref{thm:gamma} with $\Lambda = U = \R^{n}$, and hence is compact by the compactness of the graph of $K$. Since $F$ is {lower semicontinuous} {and $C_1$ and $C_2$ are closed}, the problem $(BLQOpt)$ has a solution. \qed
\end{proof}
\begin{remark}
	The hypothesis \ref{hyp2_normalized} is the normalized version of similar hypothesis used in \cite{zbMATH05351708}. Actually, it is this normalized version which should have been used in \cite{zbMATH05351708}.
\end{remark}

Now combining Theorems \ref{thm:NI} and \ref{thm:BOpt_existence}, the existence result for Single-Leader-Multi-Follower game with local responses could be proved.
\begin{proof}{\em of Theorem \ref{thm:exist_SLMFG}.}
	The continuity of $V$ follows, in view of \cite[Theorem 1.17 (c)]{zbMATH01110799}, from the continuity of $\theta_{i}$'s. Theorem \ref{thm:BOpt_existence} then yields the existence of a solution of the corresponging bilevel problem $(BLQOpt)$ with lower level represented by $f := V$. Finally, \ref{conclusion:NI-d} of Theorem \ref{thm:NI} concludes the theorem. \qed
\end{proof}

One can consider that, in Theorem \ref{thm:exist_SLMFG}, the int-dual lower semicontinuity assumption of function $F_V$ is a quite restrictive hypothesis. By the forthcoming corollary, let us show that, in some particular cases, this assumption can be naturally satisfied.

\begin{theorem}\label{thm:exist_SLMFG_linear}
	Consider a $(SLLMF)$ with $F$ being {lower semicontinuous}, $\theta_{i}$ being continuous and $\theta_{i}(\cdot,x)$ being linear for all $x \in C_{1}$ and $i = 1,\cdots,M$. Suppose that $C_{1}$ and $C_{2}$ are closed sets, that for each $i = 1,\cdots,M$, $K_{i}(\cdot,\cdot)$ is a linear translation of a fixed set such that $y\in C_{2} \mapsto K(y,x)$ is graph-compact and locally reproducible over its fixed point set for each $x \in C_{1}$, and that $FP(K(\cdot,C_{1}))$ is closed. Then the (SLLMF) admits at least a solution.
\end{theorem}

Let us clarify that, by ``$K_{i}(\cdot,\cdot)$ is a linear translation of a fixed set'': it means for each $i = 1,\cdots,M$, there exist a subset $\mathcal{K}_{i} \subset \R^{m_{i}}$ and a linear map $\mathcal{L}_{i} : \R^{m} \times \R^{n} \to \R^{m_{i}}$ such that
\begin{equation}\label{eqn:linear_translation}
	K_{i}(y_{-i},x) = \mathcal{K}_{i} + \mathcal{L}_{i}(y,x)
\end{equation}
for $(y,x) \in \R^{m} \times \R^{n}$.

\begin{proof}{\em of Theorem \ref{thm:exist_SLMFG_linear}.}
Let us first make some needful remarks about its assumption. Assume that all the assumptions of Theorem \ref{thm:exist_SLMFG_linear} hold. Then the Nikaido-Isoda function $\Psi$, as in \eqref{eqn:NI_function}, is linear for each $x$. Thus for any $x \in \R^n$ and $i = 1,\cdots,M$, there is $c_{i,x} \in \R^{m_{i}}$ in which the expression
\[
	\Psi(y,z,x) = \sum_{i=1}^{M} c_{i,x}^{T}(y_{i} - z_{i}),
\]
holds for all $y,z \in R^{m}$. Now that each $K_{i}$ is a linear translation \eqref{eqn:linear_translation}, the neighbourhood of local reproducibility can be translated uniformly over all $x \in C_{1}$ and $y \in FP(K(\cdot,x))$. This further implies that for some open set $U_{i} \subset \R^{m_{i}}$, $K(y,x) \cap U_{y,x} = \prod_{i=1}^{M} (\mathcal{K}_{i} \cap U_{i}) + \mathcal{L}_{i}(y,x)$ holds for all $x \in C_{1}$ and $y \in FP(K(\cdot,x))$. On the other hand, if $z^{\ast}(y,x)$ belongs to $\argmax_{K(y,x) \cap U_{y,x}} \Psi(y,\cdot,x)$ for each $x \in C_{1}$ and $y \in FP(K(\cdot,x))$, then there exists $\zeta^{\ast} \in \R^{m}$ for which $z^{\ast}(y,x) = \zeta^{\ast} + \mathcal{L}(y,x)$ for all such $x$ and $y$, where $\mathcal{L} := \prod_{i=1}^{M}\mathcal{L}_{i}$. Putting all of these observation together, we obtain an explicit expression for $V$ as follows
\[
	V(y,x) = \sum_{i=1}^{M} c_{i,x}^{T}(y_{i} - \zeta^{\ast} - \mathcal{L}_{i}(y,x)),
\]
which means $V$ is linear. One may observe that $\nabla_{y}V(y,x) = c_{i,x} - \hat{\mathcal{L}}_{i}^{\ast} c_{i,x}$, where $\hat{\mathcal{L}}_{i}$ comes from the coordinate decomposition $\mathcal{L}_{i}(y,x) = \hat{\mathcal{L}}_{i}y + \check{\mathcal{L}}_{i}x$. Thus the map $F_{V}(y,x) := F_{V(\cdot,x)}(y)$ reduces to the singleton $\{\nabla_{y}V(y,x)/\norm{\nabla_{y}V(y,x)}\}$ if $V(\cdot,x)$ is not a constant function and to $\{0\}$ otherwise. With the continuity of each $\theta_{i}$, one could obtain that $(F_{V},K)$ is int-dually lower semicontinuous. Now, we are in the position to prove Theorem \ref{thm:exist_SLMFG_linear}.

Consider now the problem of locally quasi-minimizing the gap function $V$ over the constraint map $K$, as in \ref{conclusion:NI-d} of Theorem \ref{thm:NI}. As described in the above discussion, all the requirements of Theorem \ref{thm:BOpt_existence} for $f(y,x) = V(y,x)$ are satisfied. Hence the graph of $x \mapsto LQOpt(V(\cdot,x),K(\cdot,x))$ is nonempty. Then by Theorem \ref{exist_quasi_opt}, problem (BLQOpt) admits at least a solution and, combining with Theorem \ref{thm:NI}, it is thus also the case for problem (SLLMF). \qed
\end{proof}

\section{Conclusions}\label{sec:Conclusion}
\DA{To conclude}, we have provided a systematic analysis regarding local solutions of Stampacchia and Minty variational and quasi-variational inequalities. Number of relations between the local and global solution sets are pointed out. We also \DA{introduced} the local concept of reproducibility for set-valued maps, which greatly expand the usability compared to its original global version. The local reproducibility allows us to solve a quasi-variational inequality locally by simplifying it into a variational inequality. Existence of local solutions of such quasi-variational inequalities were proved by combining local reproducibility and the existence for variational inequalities. The upper semicontinuity of the local solution maps of a variational inequality is also proved and used to show the upper semicontinuity of the local solution map of a quasi-optimization problem whose objective function is semistrictly quasiconvex, under the influence of local reproducibility. \blue{Pointing back to \ref{rmk:LQOpt_others} of Remark \ref{rmk:LQOpt_closedness}, it would be interesting to consider in \DA{a} forthcoming work a quantitative stability of local solutions.} Finally the stability result is applied to show existence of a solution to a Single-Leader-Multi-Follower game whose followers respond to the leader only with local generalized Nash equilibrium. This final result was approached via a lower level reformulation so that the Single-Leader-Multi-Follower game is presented as a simple bilevel problem.

\begin{acknowledgements}
The authors are grateful to the valuable comments of the reviewers which greatly improves the quality and presentation of this paper. The first author benefits from the support of PGMO \DA{through the IROE project "Numerical methods for bilevel problems: theory, numerical analysis and energy management applications (NuMeBi)"}. The second author is supported through the project ``Optimizing Leader-Follower Models in Eco-Industrial Park Management'' by King Mongkut's University of Technology Thonburi (KMUTT), Thailand Science Research and Innovation (TSRI) , and National Science, Research and Innovation Fund (NSRF) Fiscal year 2024. This research was initialized at Universit\'e de Perpignan under the support of CIMPA Research in Pairs Program during October -- December 2021.
\end{acknowledgements}



\end{document}